\documentclass{article}

\usepackage{arxiv}

\usepackage[utf8]{inputenc} 
\usepackage[T1]{fontenc}    

\usepackage{url}            
\usepackage{booktabs}       
\usepackage{amsmath, amssymb, mathtools, amsthm, amsfonts}       
\usepackage{nicefrac}       
\usepackage{microtype}      
\usepackage{graphicx}
\usepackage[table,xcdraw]{xcolor}
\definecolor{burntorange}{RGB}{191, 87, 0}
\definecolor{burntblue100}{RGB}{0, 95, 134}
\definecolor{burntgreen100}{RGB}{67, 105, 91}
\definecolor{graphred}{RGB}{160, 82, 71}
\definecolor{graphgreen}{RGB}{112, 173, 71}
\definecolor{graphblue}{RGB}{47, 82, 143}
\definecolor{tablegraybg}{RGB}{239, 239, 239}

\usepackage{tikz}
\usepackage{pgfplots}
\pgfplotsset{compat=1.16}
\tikzstyle{detector} = [rectangle, minimum width=1.0cm, minimum height=0.5cm, draw=graphred]

\usepackage{multirow}
\usepackage{natbib}
\usepackage{doi}
\usepackage{enumitem}
\usepackage{appendix}

\theoremstyle{definition}
\newtheorem{thm}{Theorem}
\newtheorem{defn}{Definition}    

\theoremstyle{remark}
\newtheorem{rem}{Remark}

\newtheorem{result}{Result}

\usepackage{dsfont}
\usepackage{bm}
\usepackage{mathrsfs}                   
\usepackage{optidef}                    

\makeatletter
\providecommand*{\diff}%
		{\@ifnextchar^{\DIfF}{\DIfF^{}}}
\def\DIfF^#1{%
		\mathop{\mathrm{\mathstrut d}}%
			\nolimits^{#1}\gobblespace}
\def\gobblespace{%
		\futurelet\diffarg\opspace}
\def\opspace{%
		\let\DiffSpace\!%
		\ifx\diffarg(%
			\let\DiffSpace\relax
		\else
			\ifx\diffarg[%
				\let\DiffSpace\relax
			\else
				\ifx\diffarg\{%
					\let\DiffSpace\relax
				\fi\fi\fi\DiffSpace}

\newcommand{\theauthors}{Sudesh K.\ Agrawal, John J.\ Hasenbein}

\newcommand{\thetitle}{Detecting Viruses in Contact Networks with Unreliable Detectors}
\newcommand{\therunningtitle}{Detecting Viruses in Contact Networks with Unreliable Detectors}

\newcommand{\thekeywords}{virus detection model; stochastic spread dynamics; unreliable detectors; Monte Carlo approximation; simulation analysis}

\title{\thetitle}

\date{} 					

\author{ \href{https://sudeshkagrawal.github.io}{Sudesh K.\ Agrawal} \\
	Graduate Program in Operations Research \& Industrial Engineering \\
	Department of Mechanical Engineering \\
	University of Texas at Austin\\
	Austin, Texas 78712 \\
	\texttt{\href{mailto:sudesh@utexas.edu}{sudesh@utexas.edu}} \\
	\And
	\href{http://sites.utexas.edu/hasenbein/}{John J.\ Hasenbein} \\
	Graduate Program in Operations Research \& Industrial Engineering \\
	Department of Mechanical Engineering \\
	University of Texas at Austin\\
	Austin, Texas 78712 \\
	\texttt{\href{mailto:jhas@mail.utexas.edu}{jhas@mail.utexas.edu}} \\
}

\renewcommand{\shorttitle}{\therunningtitle}

\usepackage{hyperref}       
\hypersetup
{
	pdfauthor={\theauthors},
	pdftitle={\thetitle},
	pdfsubject={A preprint for arXiv},
	pdfkeywords={\thekeywords},
	pdfborderstyle={/S/S/W 0},
	linkbordercolor=white,
	citebordercolor=white,
	urlbordercolor=white,
	colorlinks=false,
}

\begin{document}
\maketitle

\begin{abstract}
	This paper develops and analyzes optimization models for rapid detection of viruses in large contact networks. In the model, a virus spreads in a stochastic manner over an undirected connected graph, under various assumptions on the spread dynamics. A decision maker must place a limited number of detectors on a subset of the nodes in the graph in order to rapidly detect infection of the nodes by the virus. The objective is to determine the placement of these detectors so as to either maximize the probability of detection within a given time period or minimize the expected time to detection. Previous work in this area assumed that the detectors are perfectly reliable. In this work, it is assumed that the detectors may produce false-negative results. In computational studies, the sample average approximation method is applied to solving the problem using a mixed-integer program and a greedy heuristic. The heuristic is shown to be highly efficient and to produce high-quality solutions. In addition, it is shown that the false-negative effect can sometimes be ignored, without significant loss of solution quality, in the original optimization formulation. 
\end{abstract}


\section{Introduction}
\label{sec:intro}

In globally connected networks, computer viruses can cause significant losses in terms of both time and money. Therefore, it is essential to quickly detect them in order to implement effective mitigation measures. In this paper, we develop optimization models for virus detection in contact networks. A contact network captures information about the ``closeness'' between individuals \citep{chen2016distinction}. A virus is introduced at a node in the network through some external means and spreads because of this contact, possibly in a stochastic manner. There is some capability of detection, and we would like to know where to place a limited number of detectors to maximize the probability of detection before a given time. In this paper we consider the \emph{static detection problem} where detectors are placed in fixed locations. Related research in \citet{ding2021} and \citet{agrawal2021}, motivated by COVID-19, considers the \emph{dynamic detection problem} in which nodes are tested sequentially. While we consider specific settings in this research, because of some common mathematical features, the models we develop for the static detection problem can be applied to a more general setting of detecting an anomaly spreading through a contact network: for example, detecting contamination in a water distribution network, detecting rumors (fake news) or ads on a social network and sampling specific contents at a waste water plant.

In closely related work, Lee et al.\ (\citeyear{lee2012stochastic, lee2015optimization}) develop stochastic optimization models to study the dynamics of a virus spreading in telecommunication networks. The optimization models for rapid detection of viruses in \citet{lee2015optimization} are applicable to a broader range of spread dynamics than considered in their computational experiments. The authors consider two objectives to capture the goal of rapid detection: maximizing the probability of detection by a given time threshold, referred to as the \textbf{MPT model}, and minimizing the expected time to detection, referred to as the \textbf{MET model}. In both models, the decision maker or network manager places a limited number of perfectly reliable detectors, also called \emph{honeypots}, at time 0, before the spread starts, on a subset of nodes in the network in order to optimize one of the two objectives.

The aforementioned optimization models suffer from an exponential increase in the solution space, as the number of nodes grow large. Hence, the brute force method of evaluating every solution is computationally prohibitive for most practical networks. In \citet{lee2015optimization}, a greedy heuristic is therefore used to obtain a good solution using reasonable computational effort. Submodularity of the MPT objective function is established to ensure the greedy heuristic solution has an objective function value within a constant factor of the optimal value \citep{nemhauser1978analysis}.

This research extends the basic problem in \citet{lee2012stochastic} by considering the possibility of unreliable detectors: more specifically, detectors which produce false-negative results (but not false-positive results). We propose an MPT version for the false-negative case in Section~\ref{sec:falsenegativemodel} and focus our attention on the analysis of TN11C, RA1PC and RAEPC spread. These spread dynamics are described later using the nomenclature in Table~\ref{tab:spread_nomenclature}.
\section{Literature Review}
\label{sec:litreview}

Our model resembles that of the spread of influence through a social network, which has been studied extensively in the literature, in several domains. \citet{jackson2017influence} note that an influence model has three components: the stochastic process that describes the spread of influence, the budget constraints on creating influence, and an evaluation metric that maps the stochastic process to some real number. Domingos and Richardson (\citeyear{domingos2001mining, richardson2002mining}) propose the idea of maximizing this evaluation metric, the influence maximization problem (IMP), in the context of viral marketing. They model a market as a social network using data from ``knowledge-sharing sites,'' with customers in the market as nodes and interactions between the customers represented through Markov random fields \citep{domingos2001mining}. The model helps them determine the network value of customers and establish the effectiveness of viral marketing.

In contrast to the modeling framework of Domingos and Richardson, where the behavior of influence is modeled through a joint distribution, several operational models which describe the spread dynamics of influence also exist in the literature. These models can generally be classified as either \emph{threshold models} \citep{granovetter1978threshold}, where each node in the network has a threshold which one of its active neighbors (already-influenced neighbors) need to attain to influence it; or \emph{cascade models} \citep{goldenberg2001talk}, where each active node in the network attempts to influence its neighbors with some success probability. \citet{kempe2003maximizing} show that the influence maximization problem is NP-hard for the linear threshold model and the independent cascade model. They also outline an efficient heuristic with a theoretical bound on the optimality gap.

Solutions to the influence maximization problem in large networks often rely on the greedy algorithm proposed by \citet{nemhauser1978analysis} to efficiently find a sub-optimal solution since the influence maximization problem is NP-hard, in general. They show that the heuristic solution is within a constant factor $\frac{e-1}{e}$ of the optimal solution whenever the objective function is non-decreasing, submodular and evaluates to zero for an empty set. For solving IMP or its extensions, typically Monte-Carlo sampling is done, and a sample average approximation model is solved instead, because of an exponentially large number of possible realizations \citep{guney2019optimal}. 
\citet{guney2019efficient} proposes a linear programming relaxation with worst-case bounds that converge to the $(1-\frac{1}{e})$-bound asymptotically. 
\citet{guney2020large} model the IMP as a stochastic maximal covering location problem and use branch-and-cut of its Benders reformulation to solve the problem.

\citet{krause2008efficient} use the structure of the influence maximization problem to determine optimal strategies to efficiently deploy sensors in water distribution systems to detect the introduction of contaminants. They exploit submodularity of the problem and use a greedy heuristic to efficiently determine sensor placements on large distribution networks with theoretical performance guarantees. They consider several objectives, along with multi-criteria optimization, for which they solve for Pareto-optimal placements of sensors. They also consider two extensions: adversarial objective functions, where they monotonically transform the objective function to exploit submodularity; and robustness of sensor placements, where they model failure of sensors through a Bernoulli random variable.

In social network research several centrality measures, such as node degree, average distance and closeness, have been used when selecting nodes for spread of influence. Weighted centrality measures which linearly integrate two or more centrality measures are also used for selection of nodes \citep{fan2019information}.

Researchers have considered problems similar to the model in Lee et al.\ (\citeyear{lee2012stochastic, lee2015optimization}), with different applications.
\citet{berman1995locating} formulate a facility location problem so as to maximize the proportion of flow through the facilities, in a network with probabilistic flows. \citet{gutfraind2009optimal} develop a network interdiction model to maximize the probability of catching evaders. Both use a greedy heuristic to provide good solutions. 

This research extends the basic problem in Lee et al.\ (\citeyear{lee2012stochastic, lee2015optimization}) by considering the possibility of unreliable detectors, specifically detectors which produce false-negative results. In particular, we examine how unreliable detectors affect detection probabilities and the performance of mixed-integer programming solutions, along with greedy methods.  We are not aware of any research which incorporates the reliability of detectors, other than \citet{krause2008efficient} who consider sensor placements in water distribution networks from a robustness perspective. We modify the greedy heuristic in the literature  and perform computational studies to demonstrate the efficiency of the modified greedy heuristic, in terms of both the computational time and the solution quality.
\section{Models}
\label{sec:falsenegativemodel}

In this section we extend the model of Lee et al.\ (\citeyear{lee2012stochastic, lee2015optimization})  to model detectors which can produce false-negatives. We first introduce some notation and review the model nomenclature.

\subsection{Notation}
\label{sec:fn_notations}
Let $G(V, E)$ be a graph $G$ on a set of vertices (or nodes) $V$ with a set $E$ of \emph{undirected} edges. A virus is introduced at a node in the network by an external agent. We assume that $w_i,\, i\in V$ is the known probability mass function (pmf) of the initial location of the virus. The virus then spreads from this initial location, possibly in a stochastic manner, as described shortly. The network manager is interested in knowing where to place detectors on a subset of nodes, before the virus enters the network, in order to rapidly detect the virus. We describe the decision makers' metrics in detail in Section~\ref{sec:fn_models}.

For the virus spread dynamics, we follow the classification and terminlogy outlined in Lee et al.\ (\citeyear{lee2012stochastic, lee2015optimization}). They classify the spread dynamics using five characteristics: replication, persistence, propagation, transmissability and latency. (Transmissability refers to the probability that a node receiving the virus will be infected. This is referred to as ``Susceptibility'' in Lee et al.\ \citeyear{lee2015optimization}.) Table~\ref{tab:spread_nomenclature} reproduces the nomenclature they use to represent spread dynamics. A sequence of five characters from the table defines a model. 
The models that they develop and the extension that we develop in this section are applicable to any of the spread models listed in Table~\ref{tab:spread_nomenclature}. However, we do our computational experiments and analysis on only three of them: TN11C, RA1PC and RAEPC. TN11C is the simplest model, in which there is a virus that moves from one node to another in constant discrete time steps, but does not persist at a node that it had visited in past time steps. In contrast, in RA1PC and RAEPC, a virus replicates itself and sends out copies to other nodes, thereby infecting them, which is a more realistic scenario. While in RA1PC, a virus attempts to infect only one randomly chosen neighbor, in RAEPC it attempts to infect all its neighbors.

\begin{table}[ptbh]
    \centering
    \caption{Nomenclature for models of virus spread}
    \begin{tabular}{cl}
        \toprule
        Characteristics & \multicolumn{1}{c}{Models} \\ \midrule
        Replication & \begin{tabular}[c]{@{}l@{}}Virus transits from vertex to vertex (T),\\ Virus replicates itself and sends copies (R).\end{tabular} \\ \cmidrule(l){2-2}
        Persistence & \begin{tabular}[c]{@{}l@{}}Only newly infected vertices distribute the virus (N),\\ All infected vertices distribute the virus (A).\end{tabular} \\ \cmidrule(l){2-2}
        Propagation & \begin{tabular}[c]{@{}l@{}}Virus propagates to \\ one randomly selected neighbor (1),\\ every neighbor (E).\end{tabular} \\ \cmidrule(l){2-2}
        Transmissability & \begin{tabular}[c]{@{}l@{}}Virus is transmissable with probability 1 (1),\\ Virus is transmissable with probability $p<1$ (P).\end{tabular} \\ \cmidrule(l){2-2}
        Latency & \begin{tabular}[c]{@{}l@{}}Transmission occurs\\ in constant unit time steps (C),\\ according to an exponential distribution (M),\\ according to a general distribution (G).\end{tabular} \\ \cmidrule(l){2-2}
    \end{tabular}
    \label{tab:spread_nomenclature}
\end{table}

We now introduce two definitions that are needed in the sequel. 
\begin{defn}[Submodular function]
Let $2^V$ denote the power set of a set $V$. A function $f: 2^V \rightarrow \mathds{R}$ is a \emph{submodular function} if $f(S_1 \cup \{j\}) - f(S_1) \geq f(S_2 \cup \{j\}) - f(S_2)\quad\forall\,\,\, S_1 \subseteq S_2 \subseteq V,\,\,\, j\in V\setminus S_2$.
\end{defn}
\noindent Roughly speaking, a set function is submodular if it has a diminishing marginal utility in the size of the set. 
\begin{defn}[Supermodular function]
Let $2^V$ denote the power set of a set $V$. A function $f: 2^V \rightarrow \mathds{R}$ is a \emph{supermodular function} if $-f(\cdot)$ is submodular.
\end{defn}

\subsection{Optimization formulations}
\label{sec:fn_models}

In this section we discuss the MPT formulation and the MET formulation developed in \cite{lee2012stochastic}, and then extend the MPT formulation for the case where detectors with false-negative results are allowed.

In the MPT formulation the objective is to choose a set of detectors $S \subseteq V$ to maximize the probability of detecting the virus by a given time threshold $t_0$, with $|S|=k$. We assume $|V|\geq k$. Let $T_S$ be the first-passage time for the virus to hit the set $S$. The MPT optimization formulation is:
\begin{maxi!}|s|[2]    
{\substack{S\subseteq V}}
{f(S) = \mathds{P}(T_S \leq t_0)}  
{\label{eq:mptmodel}} 
{z_P^\ast =}  
\addConstraint{|S|}{=k.}
\end{maxi!}

For any of the spread models we consider, the virus spread dynamics form a stochastic process. A \emph{realization} of this process is the set of infected nodes up to some time $t_0$.
Notice that different sample paths of virus spreads may correspond to the same realization, since we do not need to record the order in which nodes are infected. 
If we enumerate all the realizations, then we can reformulate (\ref{eq:mptmodel}) as an integer program (IP):

\begin{maxi!}|s|[2]    
{\substack{x,\,y}}
{\sum\limits_{\omega\in\Omega} \phi_\omega y_\omega}  
{\label{eq:mptmodelreform}} 
{z_P^\ast =}  
\addConstraint{\sum\limits_{j\in V} x_j}{=k\label{eq:mptmodelreform_const1}}
\addConstraint{y_\omega}{\leq \sum\limits_{j\in P_\omega} x_j,\label{eq:mptmodelreform_const2}}{\quad \omega\in\Omega}
\addConstraint{x_j}{\in \{0, 1\},\label{eq:mptmodelreform_const3}}{\quad j\in V}
\addConstraint{y_\omega}{\in \{0, 1\},\label{eq:mptmodelreform_const4}}{\quad \omega\in\Omega,}
\end{maxi!}
where
\begin{itemize}[leftmargin=72pt]
    \item[$\Omega :$] set of indices $\omega$ of all realizations of the virus spread
    \item[$P_\omega \subseteq V:$] 
    set of infected vertices in realization $\omega$
    \item[$\phi_\omega :$] probability of realization $\omega$ ($\phi$ is the p.m.f.\ over $\omega$.)
    \item[$y_\omega,\ \omega\in\Omega:$] binary variables
    \item[$x_j,\ j\in V:$] 1 if node $j$ is selected as a detector, 0 otherwise.
\end{itemize}

Constraint (\ref{eq:mptmodelreform_const1}) limits the number of available detectors to $k$, and constraint (\ref{eq:mptmodelreform_const2}) ensures that $y_\omega = 1$ only if the realization indexed by $\omega$ contains at least one honyepot. Given (\ref{eq:mptmodelreform_const3}), we can relax (\ref{eq:mptmodelreform_const4}) to
$y_{\omega} \in [0,1]$, yielding a mixed integer program (MIP). In our implementation we solve the problem as an IP since we found the solver took slightly less time when solved as an IP.

In the MET formulation the objective is to choose a set of detectors $S \subseteq V$, with $|S|=k$, to minimize the expected time until we detect the virus. The MET optimization formulation is:
\begin{mini!}|s|[2]    
{\substack{S\subseteq V}}
{f(S) = \mathds{E}T_S}  
{\label{eq:metmodel}} 
{z_E^\ast =}  
\addConstraint{|S|}{=k.}
\end{mini!}

In \cite{lee2012stochastic}, it is observed that (\ref{eq:mptmodelreform}) and (\ref{eq:metmodel}) apply to all of the spread models introduced earlier since, for purposes of the MPT and the MET, each model is defined only by the set of realizations $P_\omega$ and the pmf on this set. We reproduce their two main results in Theorem~\ref{thm:mptsubmodular} and Theorem~\ref{thm:metsupermodular}.

\begin{thm}\label{thm:mptsubmodular}(\cite{lee2015optimization})
Let $G=(V, E)$ be connected with $|V| < \infty$. Let $S \subset V$ be a set of detectors, let $T_S$ be the hitting time of the stochastic process governing virus spread to set $S$, and let $f(S) = \mathds{P}(T_S \leq t_0)$ denote the probability of detecting the virus by time $t_0 > 0$. Then, $f$ is submodular.
\end{thm}

\begin{thm}\label{thm:metsupermodular}(\cite{lee2015optimization})
Under the hypotheses of Theorem~\ref{thm:mptsubmodular}, let $g(S)=\mathds{E}T_S$ denote the expected time to detect the virus, and assume $g(S)<\infty$. Then, $g$ is supermodular.
\end{thm}

While the MPT and the MET formulations can be used for any of the spread dynamics discussed, enumerating all realizations is computationally intractable for big networks because $|\Omega|$ is exponentially large in the number of nodes in the graph. The MIP reformulation lends itself to a \emph{sample average approximation} (SAA) formulation, which is more tractable. We can simulate $n$ sample paths of the virus spread instead of enumerating all realizations and solve (\ref{eq:mptmodelreform}) on these samples, leading to the following sample average approximation formulation (\ref{eq:mptmodelreform_saa}):

\begin{itemize}[leftmargin=0in]
\item[] (\textbf{SAA formulation for completely reliable detectors})
\begin{maxi!}|s|[2]    
    {\substack{x,\,y}}
    {\frac{1}{n}\sum\limits_{l=1}^{n} y_l}  
    {\label{eq:mptmodelreform_saa}} 
        {z_n^\ast =}  
    \addConstraint{\sum\limits_{j\in V} x_j}{=k\label{eq:mptmodelreform_saa_const1}}
    \addConstraint{y_l}{\leq \sum\limits_{j\in P_l} x_j,\label{eq:mptmodelreform_saa_const2}}{\quad l=1,2,\ldots n}
    \addConstraint{x_j}{\in \{0, 1\},\label{eq:mptmodelreform_saa_const3}}{\quad j\in V}
    \addConstraint{y_l}{\in \{0, 1\},\label{eq:mptmodelreform_saa_const4}}{\quad l=1,2,\ldots n,}
\end{maxi!}
\end{itemize}
where $l$ indexes the $n$ sample paths. 
Under some technical conditions, it is well-known that SAA solutions converge to the solutions value, i.e., $z_n \longrightarrow z_p$. Furthermore, it can be shown that $\mathbb{E}z_n^\ast - z_P^\ast \geq 0, \,\forall n$ (positive bias); and $\mathbb{E}z_n^\ast \geq \mathbb{E}z_{n+1}^\ast, \,\forall n$ i.e., the bias is smaller for larger samples (see \cite{mak1999monte, norkin1998branch} for details).

A greedy heuristic is described in \cite{lee2015optimization}. Submodularity of the objective function ensures that the greedy heuristic produces a solution whose value is within a constant factor $\frac{e-1}{e}$ of the optimal objective value \citep{nemhauser1978analysis}.

 We now develop an MPT formulation for the false-negative case, where the detector at a given node fails to detect the virus present on that node with probability $r_j$ with $r_j \in [0, 1),\forall j\in V$. To describe the formulation and simulation methods for models with detectors that can yield false negatives, we now introduce some definitions. 

\begin{defn}[Virtual detection]
    Let $\tilde{R}_j$ be a Bernoulli random variable associated with a node $j$, with probability of success $1-r_j$. Let $\tilde{R}_j=1$ when the success event occurs, and $0$ otherwise. A \emph{virtual detection} is said to occur at node $j$ at a given time when $\tilde{R}_j=1$ at that time.
\end{defn}
\begin{defn}[Successful detection]
    A \emph{successful detection} is said to occur at a node with a detector at a given time, when there is a virus at that node at that time, and the detector signals the presence of that virus.
\end{defn}

Note that virtual detection is not concerned with the presence of a detector; it just indicates that if a detector is present at a node, then a successful detection occurs at that node. These definitions will be helpful in writing the MPT formulation as an integer program and describing our modified  greedy heuristic.

\begin{rem}
In our implementation of the MPT model, we allow multiple independent coin flips at a node (for virtual detection) in the TN11C spread model if the virus revisits the node, but to limit book-keeping we do a coin flip only once at the infected nodes in the RA1PC and the RAEPC spread models. For these latter two models, this implies that a detector only has one ``chance'' to detect a virus present at a node.
\end{rem}

When detectors are not completely reliable, the virus may
not trigger an alert at a node with a detector installed. 
 So, we need to redefine our optimization problem. Let $T^d_S$ be the first-detection time after the virus hits the set $S$. (Note that $T^d_S=T_S$ when detectors are completely reliable.) The false-negative version of MPT can be formulated as:
\begin{maxi!}|s|[2]    
{\substack{S\subseteq V}}
{f(S) = \mathds{P}(T^d_S \leq t_0)}  
{\label{eq:mptfnmodel}} 
{z_{P_{\text{FN}}}^\ast =}  
\addConstraint{|S|}{=k,}
\end{maxi!}
where the subscript `FN' is used to indicate the false-negative case.

\begin{thm}\label{thm:fnmptsubmodular}
Let $G=(V, E)$ be connected with $|V| < \infty$. Let $S \subset V$ be a set of detector nodes, and $T^d_S$ be the first-detection time using set $S$. Furthermore, let $f(S) = \mathds{P}(T^d_S \leq t_0)$ denote the probability of successfully detecting the virus by time $t_0 > 0$. Then, $f(\cdot)$ is submodular.
\end{thm}
\proof
With the revised definition of $f(S)$, note that $f(S\cup \{j\}) - f(S)$ sums values of $\phi_\omega$ for the realizations in which the detector at $j$ successfully detected the virus, but in which none of the detectors in $S$ successfully detected the virus. Next, given $S_1 \subseteq S_2$ and $j \notin S_2$ we then have 
$$f(S_1 \cup \{j\}) - f(S_1)  \geq f(S_2 \cup \{j\}) - f(S_2).$$
The inequality clearly holds if $j \in S_1$ or $j \in S_2$. Suppose then that $j \in S_2^C$. Then, the sum defining right-hand side consists of realizations in which none of the detectors in
$S_2$ (or $S_1$) successfully detected the virus, but the detector at $j$ did. Notice then that all such realizations 
are also contained in the sum defining the left-hand side, since
if no detector in $S_2$ detected the virus, neither did any detector in $S_1$. Thus, the left-hand side contains at least the same summands that the right-hand side does.
\endproof

\begin{rem}
    The proof holds for any of the spread dynamics outlined in Table~\ref{tab:spread_nomenclature}.
\end{rem}

\noindent Next, we present a similar result pertaining to the MET formulation, with false negatives. 

\begin{thm}\label{thm:fnmetsupermodular}
Under the hypothesis of Theorem~\ref{thm:fnmptsubmodular}, let $g(S)=\mathds{E}T^d_S$ denote the expected time to detect the virus, and assume $g(S)<\infty$. Then, $g(\cdot)$ is supermodular.
\end{thm}
\proof{Proof:}
Let $F(S, t) = \mathds{P}(T^d_S \leq t)$ denote the probability of detecting the virus by installing the detectors at nodes in the set $S$, by the general threshold time $t$. Then,
\begin{equation}
g(S) = \mathds{E}T^d_S = \int\limits_0^\infty\left [ 1 - F(S, t) \right]\diff t.
\end{equation}
%
From Theorem~\ref{thm:fnmptsubmodular} we know that $F(\cdot, t)$ is submodular for any $t>0$, i.e.,
\begin{equation}\label{eq:submodularitydefintheorem2}
F(S_1 \cup \{j\}, t) - F(S_1, t) \geq F(S_2 \cup \{j\}, t) - F(S_2, t)\quad\forall\,\,\, S_1 \subseteq S_2 \subseteq V,\,\,\, j\in V\setminus S_2.
\end{equation}
So, 
\begin{equation*}
    \begin{array}{rcl}
g(S_1 \cup \{j\})-g(S_1) & = & \int\limits_0^\infty [1-F(S_1 \cup \{j\}, t)]\diff t - \int\limits_0^\infty [1-F(S_1, t)]\diff t \\
     & = & \int\limits_0^\infty [F(S_1, t)-F(S_1 \cup \{j\}, t)]\diff t \\
     & \leq & \int\limits_0^\infty [F(S_2, t)-F(S_2 \cup \{j\}, t)]\diff t \\
         & = & g(S_2 \cup \{j\}) - g(S_2), 
    \end{array}
\end{equation*}
\noindent where the inequality follows from~(\ref{eq:submodularitydefintheorem2}).
\endproof

\begin{rem} We can use a backward greedy algorithm for a constant-factor guarantee in the MET problem. However, we restrict our focus in this paper on the MPT version and the corresponding standard greedy heuristic.
\end{rem}

Now, returning to the MPT problem, note that we can 
can reformulate (\ref{eq:mptfnmodel}) as a mixed-intger program:
\begin{maxi!}|s|[2]    
{\substack{x,\,u}}
{\sum\limits_{\omega\in\Omega} \phi_\omega u_\omega}  
{\label{eq:mptfnmodelreform}} 
{z_{P_{\text{FN}}}^\ast =}  
\addConstraint{\sum\limits_{j\in V} x_j}{=k\label{eq:mptfnmodelreform_const1}}
\addConstraint{u_\omega}{\leq \sum\limits_{j\in P_{\text{FN},\,\,\omega}} x_j,\label{eq:mptfnmodelreform_const2}}{\quad \omega\in\Omega}
\addConstraint{x_j}{\in \{0, 1\},\label{eq:mptfnmodelreform_const3}}{\quad j\in V}
\addConstraint{u_\omega}{\in \{0, 1\},\label{eq:mptfnmodelreform_const4}}{\quad \omega\in\Omega,}
\end{maxi!}
where
\begin{itemize}[leftmargin=72pt]
    \item[$\Omega :$] set of indices $\omega$ of all realizations of the virus spread
    \item[$P_{\text{FN},\,\,\omega} \subseteq V:$] set of vertices that virtually detect the virus up to $t_0$ in realization $\omega$
    \item[$\phi_\omega :$] probability of realization $\omega$ occurring ($\phi$ is the p.m.f.\ over $\omega$.)
    \item[$u_\omega,\  \omega\in\Omega:$] binary variables
    \item[$x_j,\ j\in V:$] binary variables to determine the set of detector nodes.
\end{itemize}
As before, we can formulate a SAA version of the problem, which is used in the next section: 
\begin{itemize}[leftmargin=0in]
\item[] (\textbf{SAA formulation for detectors with false negatives})
\begin{maxi!}|s|[2]    
    {\substack{x,\, u}}
    {\frac{1}{n}\sum\limits_{l=1}^{n} u_l}  
    {\label{eq:mptfnmodelreform_saa}} 
    {z_{\text{FN},\,\,n}^\ast =}  
    \addConstraint{\sum\limits_{j\in V} x_j}{=k\label{eq:mptfnmodelreform_saa_const1}}
    \addConstraint{u_l}{\leq \sum\limits_{j\in P_{\text{FN},\,\,l}} x_j,\label{eq:mptfnmodelreform_saa_const2}}{\quad l=1,2,\ldots n}
    \addConstraint{x_j}{\in \{0, 1\},\label{eq:mptfnmodelreform_saa_const3}}{\quad j\in V}
    \addConstraint{u_l}{\in \{0, 1\},\label{eq:mptfnmodelreform_saa_const4}}{\quad l=1,2,\ldots n.}
\end{maxi!}
\end{itemize}

\section{Simulation and Greedy Heuristic}
\label{sec:fn_simandheuristic}
In our computational studies in Section~\ref{sec:fn_computationalstudies} we focus our analyses on TN11C, RA1PC and RAEPC spread dynamics. In this section we describe how to generate Monte Carlo samples by simulating the stochastic spread given these virus dynamics. We also describe the greedy heuristic which is used in conjunction with the SAA formulations (\ref{eq:mptmodelreform_saa}) and (\ref{eq:mptfnmodelreform_saa}). The following inputs are required for any of the stochastic spread models:
\begin{itemize}[leftmargin=0.25in]
    \item the time threshold $t_0$,
    \item the number of sample paths $n$ to be generated,
    \item and for the false negative case, the false negative probability $r$. (For simplicity, we assume $r_j = r, \,\forall j\in V$.)
\end{itemize}

Recall that in the TN11C model, the virus hops from one node to another in constant time steps. Let $L$ be the set of sampled virus paths. For each sample path $l \in L$, at $t = 0$, we first choose an initial node for the virus to infect, uniformly at random among all nodes in $V$. At $t=1, 2, \ldots t_0$, we uniformly randomly pick a neighbor of the currently infected node for the virus to infect. We keep track of the infected nodes and generate a \emph{virus spread matrix} (VSM) in which the rows represent sample paths, and the columns correspond to time. Hence, the VSM has $|L|$ rows and $t_0+1$ columns. A sample VSM for $|L|=3$ and $t_0=2$ for a TN11C model is as follows:
\begin{equation}\label{eq:fn_tn11c_vsm_example}
    \begin{bmatrix}
    5 & 192 & 3 \\
    4 & 3 & 5 \\
    1 & 11 & 13 \\
    \end{bmatrix}.
\end{equation}

When detectors can produce false-negative results, we additionally generate a \emph{virtual detection matrix} (VDM) whose elements take values in $\{0, 1\}$. An element of this matrix indicates that a detector placed at the node corresponding to that element would detect the virus (i.e., a true positive would occur). The VDM also has $|L|$ rows and $t_0+1$ columns. To generate this matrix, we simulate independent Bernoulli random variables $\tilde{R}_j$ for each node and for each point of time in the VSM. Note that this means we check for virtual detection every time the virus visits a node. For example, if the virus visits node $i$ twice in a sample path, then we generate two independent samples of $\tilde{R}_{i}$. An example VDM for the VSM example in (\ref{eq:fn_tn11c_vsm_example}) appears below:
\begin{equation*}
    \begin{bmatrix}
    1 & 1 & 1 \\
    1 & 0 & 1 \\
    0 & 0 & 1 \\
    \end{bmatrix}.
\end{equation*}

In the RA1PC model, the virus at an infected node tries to replicate and send a copy of itself to a randomly chosen neighbor, at every time step. For each sample path $l \in L$, at $t = 0$, we again choose an initial node for the virus to infect uniformly at random. At $t=1, 2, \ldots t_0$, we uniformly randomly pick a neighbor of each of the currently infected nodes for the virus to infect, with successful infection happening independently with probability $p$. We keep track of the infected nodes and generate a list of $|L|$ rows, where each row is a set of infected nodes. The rows in this case need not contain the same number of elements. Hence, the resulting matrix is, in general, what is sometimes referred to as a ``ragged matrix.'' As such, we still refer to this object as a VSM. An example VSM for $|L|=3$ and $t_0=2$ is:
\begin{equation}\label{eq:fn_ra1pc_vsm_example}
    \begin{bmatrix}
    5 & 192 & 3 & 7 & 13 \\
    4 & 3 & 5 \\
    1 & 11 & 13 & 23 \\
    \end{bmatrix}.
\end{equation}

When detectors can produce false-negative results, we additionally generate a VDM by simulating independent Bernoulli random variables $\tilde{R}_j$ for each infected node in each of the $|L|$ rows of the VSM. Note that this means we check for virtual detection only once per node since we store the infected nodes in a sample path as a set. An example VDM for the VSM example in (\ref{eq:fn_ra1pc_vsm_example}) is:
\begin{equation}\label{eq:fn_ra1pc_vdm_example}
    \begin{bmatrix}
    1 & 0 & 1 & 0 & 1\\
    0 & 1 & 1 \\
    0 & 1 & 1 & 1\\
    \end{bmatrix}.
\end{equation}

The RAEPC simulation process is analogous to the RA1PC except
for the spread dynamics. In RAEPC, the virus at an infected node tries to replicate and send a copy of itself to all of its neighbors, at every time step. 
As before, we keep track of the infected nodes and generate a list of $|L|$ rows, where each element is a set of infected nodes. 
When detectors can produce false-negative results, we additionally generate a VDM by simulating a Bernoulli random variable $\tilde{R}_j$ for each infected node in each of the $|L|$ rows. 



\noindent\textbf{Greedy heuristic:} 
We assume now that there are a finite number of nodes that are labeled 1, 2, \ldots, $|V|$.
A Hadamard product (element-wise multiplication) of the VSM and the VDM gives us another matrix: the \emph{successful detection matrix} (SDM). The SDM corresponding to the VSM in (\ref{eq:fn_ra1pc_vsm_example}) and the VDM in (\ref{eq:fn_ra1pc_vdm_example}) is as follows:
\begin{equation*}\label{eq:fn_ra1pc_sdm_example}
    \begin{bmatrix}
    5 & 0 & 3 & 0 & 13 \\
    0 & 3 & 5 \\
    0 & 11 & 13 & 23 \\
    \end{bmatrix}.
\end{equation*}
When the detectors are completely reliable i.e., $r=0$, then the VDM contains only 1's, and the SDM is equal to the VSM. Construction of the SDM is useful in describing the greedy heuristic. Each row in an SDM represents the set of nodes in the corresponding sample path where a successful detection could occur, if a detector is indeed placed at the corresponding node. 

The greedy algorithm for selecting detector nodes for the MPT objective, with false negatives works as follows. For any spread model, generate the VSM and VDM with $|L|$ rows, corresponding to the number of samples paths in the simulation. The SDM is then formed as described above. We initialize $S_0$ to be the empty set. The first detector selected is a node $j_0$ which is present in the largest number of rows in the SDM, with ties broken arbitrarily. Set $S_1 = S_0 \cup \{j_0\}.$ Now at iteration $i$, a row is a \emph{candidate row} if none of the nodes in the row
is contained in $S_i$. Then in iteration $i$ we select a node $j_i$ which is present in the largest number of candidate rows, with ties broken arbitrarily, and set $S_{i+1}= S_i \cup \{j_i\}.$ The algorithm terminates when $S_k$ is determined. This set $S_k$ is then set to be $S$ for the MPT objective.




This greedy heuristic may be preferred over the MIP with SAA if it can give us an ``acceptable'' solution more quickly, especially for large networks. In Section~\ref{sec:fn_computationalstudies} we undertake computational studies to analyze how the heuristic fares against a solution via the MIP formulation, both in terms of computational time and solution quality.

\section{Computational Studies}
\label{sec:fn_computationalstudies}

\citet{lee2012stochastic} uses real data from an Asian telecommunication company to build contact networks and uses a down-sampling scheme, $c$-core decomposition \citep{seidman1983network}, to build network samples of different scales. He also discusses how the decomposition scheme keeps the relevant properties of the network intact. 
A $c$-core of a connected graph is not necessarily connected. Lee considers only the largest connected component of the $c$-core decomposed network, because if we do not install a detector in an isolated component, then we cannot detect a virus in that component. 

We use a network generated using email data from a research institution \citep{snapnets, leskovec2007graph}, and therefore the context of our computational studies is the spread of email viruses. This network could be viewed as representing a collaboration network of researchers, which has been shown to behave like a ``small world'' network \citep{newman2001structure}. In a small world network the typical path length $M$ between two randomly chosen nodes is small compared to the number of nodes $|V|$ in the graph: more precisely, $M=\mathcal{O}(\log |V|)$ \citep{watts1998collective}.
We use the largest connected component of the 6-core of the aforementioned network, removing all self-loops, and refer to this decomposed network as the ``6-core EU email network,'' or just the ``EU email network'' when it is clear that we mean the 6-core version of the network. This network has 5400 nodes and 73315 edges with an average degree of about 27. Details of the network can be found in Appendix~\ref{sec:euemailnetwork_appendix}.

We use simulation to determine a suitable value for $t_0$ to ensure quick detection of a virus. A reasonable value depends on how quickly the virus spreads over the network and on the decision maker's tolerance. In our computational studies, we use $t_0=4$ for TN1PC, $t_0=3$ for RA1PC and $t_0=1$ for RAEPC, corresponding to infection of approximately 0.1\% of the nodes when $p=1$. Details of the simulation used to determine these values are in Appendix~\ref{sec:fn_choosingtime_appendix}. We perform the following computational studies:
\begin{enumerate}[leftmargin=0.25in]
    \item comparing the solutions from MIP and greedy heuristic
    \item evaluating the quality of the greedy solution using the multiple replication procedure (MRP)
    \item studying the effect of the false-negative rate on the probability of detection
    \item studying the effect of ignoring detector fallibility.
\end{enumerate}

All the experiments were done on a 64-bit system with Intel(R) Core(TM) i7-10875H CPU @ 2.30GHz 2.30 GHz with 31.8 GB usable RAM and 8 cores (16 logical processors).

\subsection{Comparing MIP and Greedy Solutions}
\label{sec:fn_compstudy_comparemipandgreedy}

For this study we generate $n$ samples for the RA11C spread model and solve (\ref{eq:mptfnmodelreform_saa}) using these samples. We use the Gurobi solver \citep{gurobi} via its Java API to solve the MIP and terminate when it achieves a relative gap of 1\% between the best integer objective $\underline{z}^{\ast}_n$ and the best linear-programming-relaxed objective $\overline{z}^{\ast}_n$, or when the time expended exceeds 1 hour. The relative gap provides a bound on the accuracy of the MIP optimal solution. The integer feasibility tolerance is set to $1.0\cdot10^{-9}$, and only one processor thread is used. 
%
We refer to the solution $S^\ast_n$ produced by the solver as ``MIP.'' We also use the greedy heuristic on the same samples and obtain the objective value $z^h_n$. We refer to the solution $S^h_n$ from the greedy heuristic as ``Greedy.'' We note the computational times for both.
%

We also examine the difference $\mu = f(S^h_n) - f(S^\ast_n)$ in the objective value of Greedy and MIP and form the corresponding confidence interval using McNemar's procedure \citep{mcnemar1947note}.
A (1-$\alpha$)-level CI on the difference between marginal proportions is given by:
\begin{equation}\label{eq:mcnemar_marginalpropdiff}
    \hat{d} = \frac{n_{12}}{n^{\prime\prime}} - \frac{n_{21}}{n^{\prime\prime}},
\end{equation}
with sample variance:
\begin{equation}
    s^2(\hat{d}) = \frac{1}{n^{\prime\prime}}\cdot \left [\frac{n_{12}}{n^{\prime\prime}}(1 - \frac{n_{12}}{n^{\prime\prime}}) + \frac{n_{21}}{n^{\prime\prime}}(1 - \frac{n_{21}}{n^{\prime\prime}}) + 2\frac{n_{12}n_{21}}{(n^{\prime\prime})^2} \right ],
\end{equation}
where $n_{11}$, $n_{12}$, $n_{21}$, and $n_{22}$ are described in Table~\ref{tab:mcnemars_table}, and $n^{\prime\prime}$ is the total number of sample paths generated in evaluating the greedy and the MIP solution. The approximate (1-$\alpha$)-level CI is given by
\begin{equation}
    \left[ \hat{d} - z_{\frac{\alpha}{2}}s(\hat{d}), \hat{d} + z_{\frac{\alpha}{2}}s(\hat{d}) \right],
\end{equation}
where $z_\alpha$ is the $\alpha$ quantile of the standard normal. In (\ref{eq:mcnemar_marginalpropdiff}) we subtract the ``MIP positive'' proportion from the ``Greedy positive'' proportion, which means if a CI covers strictly positive values, then we infer that there is strong evidence that the greedy solution is better. If a CI covers strictly negative values, then there is evidence that the MIP solution is better. 
\begin{table}[ptbh]
    \centering
    \caption{``Greedy positive'' indicates that with the greedy solution $S^h_n$, the sample realization is detected, and ``Greedy negative'' indicates the virus is undetected with the greedy solution. Similarly, ``MIP positive'' and ``MIP negative'' indicate whether or not the virus is detected in the MIP solution.}
    \begin{tabular}{cccc}
        \toprule
         & \textbf{MIP positive} & \textbf{MIP negative} & \textbf{Total} \\ \midrule
        \textbf{Greedy positive} & $n_{11}$ & $n_{12}$ & $n_{11}+n_{12}$ \\
        \textbf{Greedy negative} & $n_{21}$ & $n_{22}$ & $n_{21}+n_{22}$ \\
        \textbf{Total} & $n_{11}+n_{21}$ & $n_{12}+n_{22}$ & $n^{\prime\prime}$ \\ 
        \bottomrule
    \end{tabular}
    \label{tab:mcnemars_table}
\end{table}

The results appear in Table~\ref{tab:ra11c_compare_mip_and_greedy_t3_r5}. In the table, if $z^h_n > \underline{z}^{\ast}_n$, then the heuristic provides a better solution on the relatively  small set of sampled paths indicated in the second column of the table. Otherwise, the MIP solution is better. Highlighted cells in the table indicate cases in which the solver terminated due to the time limit of 1 hour. We generate CIs on the difference $\mu = f(S^h_n) - f(S^\ast_n)$ corresponding to 95\% significance, using 2 million samples ($n^{\prime\prime}=2000000$). 

\begin{rem}
Unless stated otherwise, all further analyses use MIP solutions attainable within 1 hour or approximately 1\% relative gap.
\end{rem}

\begin{table}[ptbh]
    \centering
    \caption{Comparison of MIP and Greedy, under RA11C spread on 6-core EU email network with $t_0=3$ and $r=0.05$: computational time and 95\% confidence interval on the difference $\mu = f(S^h_n) - f(S^\ast_n)$ using McNemar's procedure ($n^{\prime\prime}=2000000$). A positive difference indicates that the greedy solution performed better. Note that columns $\bm{\underline{z}^{\ast}_n}$, $\bm{\overline{z}^{\ast}_n}$, $\bm{z^h_n}$ and both wall times are not from the McNemar's procedure.}
    \label{tab:ra11c_compare_mip_and_greedy_t3_r5}
    \begin{tabular}{@{}cccccccc@{}}
        \toprule
        \multicolumn{2}{c}{\textbf{}} & \multicolumn{3}{|c|}{\textbf{MIP}} & \multicolumn{2}{c|}{\textbf{Greedy}} & $\bm{f(S^h_n) - f(S^\ast_n)}$ \\ \midrule
        $\bm{k}$ & $\bm{n}$ & $\bm{\underline{z}^{\ast}_n}$ & $\bm{\overline{z}^{\ast}_n}$ & \textbf{\begin{tabular}[c]{@{}c@{}}Wall time\\ (sec)\end{tabular}} & $\bm{z^h_n}$ & \textbf{\begin{tabular}[c]{@{}c@{}}Wall time\\ (sec)\end{tabular}} & \textbf{95\% CI} \\ \midrule
         & 1000 & 0.5990 & 0.6040 & 0.3901 & 0.5970 & 0.0892 & {[}0.0258, 0.0267{]} \\
         & 5000 & 0.5350 & 0.5402 & 7.434 & 0.5374 & 0.1022 & {[}-0.0040, -0.0033{]} \\
         & 10000 & 0.5263 & 0.5312 & 27.28 & 0.5277 & 0.1526 & {[}0.0019, 0.0023{]} \\
         & 30000 & \cellcolor[HTML]{EFEFEF}0.5249 & \cellcolor[HTML]{EFEFEF}0.5305 & \cellcolor[HTML]{EFEFEF}3600 & 0.5250 & 0.5054 & \colorbox{tablegraybg}{{[}0.0005, 0.0008{]}} \\
        \multirow{-5}{*}{50} & 50000 & \cellcolor[HTML]{EFEFEF}0.5201 & \cellcolor[HTML]{EFEFEF}0.5316 & \cellcolor[HTML]{EFEFEF}3600 & 0.5239 & 0.8164 & \colorbox{tablegraybg}{{[}0.0031, 0.0039{]}} \\ \cmidrule(l){2-8}
         & 1000 & 0.8420 & 0.8500 & 3.754 & 0.8330 & 0.1300 & {[}0.0173, 0.0184{]} \\
         & 5000 & \cellcolor[HTML]{EFEFEF}0.7506 & \cellcolor[HTML]{EFEFEF}0.7624 & \cellcolor[HTML]{EFEFEF}3600 & 0.7484 & 0.1448 & \colorbox{tablegraybg}{{[}0.0031, 0.0039{]}} \\
         & 10000 & \cellcolor[HTML]{EFEFEF}0.7345 & \cellcolor[HTML]{EFEFEF}0.7501 & \cellcolor[HTML]{EFEFEF}3600 & 0.7350 & 0.2722 & \colorbox{tablegraybg}{{[}-0.0021, -0.0016{]}} \\
         & 30000 & \cellcolor[HTML]{EFEFEF}0.7247 & \cellcolor[HTML]{EFEFEF}0.7560 & \cellcolor[HTML]{EFEFEF}3600 & 0.7303 & 0.9234 & \colorbox{tablegraybg}{{[}0.0051, 0.0058{]}} \\
        \multirow{-5}{*}{100} & 50000 & \cellcolor[HTML]{EFEFEF}0.7092 & \cellcolor[HTML]{EFEFEF}0.7578 & \cellcolor[HTML]{EFEFEF}3600 & 0.7301 & 1.4062 & \colorbox{tablegraybg}{{[}0.0188, 0.0198{]}} \\ \cmidrule(l){2-8}
         & 1000 & 1 & 1 & 0.2053 & 1 & 0.1570 & {[}0.0522, 0.0535{]} \\
         & 5000 & \cellcolor[HTML]{EFEFEF}0.9290 & \cellcolor[HTML]{EFEFEF}0.9512 & \cellcolor[HTML]{EFEFEF}3600 & 0.9270 & 0.2304 & \colorbox{tablegraybg}{{[}-0.0002, 0.0006{]}} \\
         & 10000 & \cellcolor[HTML]{EFEFEF}0.9084 & \cellcolor[HTML]{EFEFEF}0.9369 & \cellcolor[HTML]{EFEFEF}3600 & 0.9098 & 0.4400 & \colorbox{tablegraybg}{{[}0.0004, 0.0011{]}} \\
         & 30000 & \cellcolor[HTML]{EFEFEF}0.8717 & \cellcolor[HTML]{EFEFEF}0.9306 & \cellcolor[HTML]{EFEFEF}3600 & 0.9007 & 1.268 & \colorbox{tablegraybg}{{[}0.0251, 0.0261{]}} \\
        \multirow{-5}{*}{200} & 50000 & \cellcolor[HTML]{EFEFEF}0.8675 & \cellcolor[HTML]{EFEFEF}0.9306 & \cellcolor[HTML]{EFEFEF}3600 & 0.8973 & 2.051 & \colorbox{tablegraybg}{{[}0.0294, 0.0303{]}} \\ \cmidrule(l){2-8}
    \end{tabular}
\end{table}

\noindent From Table~\ref{tab:ra11c_compare_mip_and_greedy_t3_r5} we make the following observations:
\begin{enumerate}[leftmargin=0.25in]
    \item The greedy heuristic outperforms the MIP in terms of computational time for larger sample sizes $n$. 
    \item For larger values of $n$ or $k$ the solver is not able to solve the MIP to optimality (within the given tolerance) within an hour.
    \item There is a diminishing marginal value of adding detectors.
    \item For the majority of cases, the confidence interval covers strictly positive values, indicating that the greedy heuristic is likely better most of the time. Even in cases where the CI does not cover strictly positive values, we see that the MIP and greedy solution are the same in the first two significant digits. 
\end{enumerate}
Overall, this provides compelling evidence that the greedy heuristic can provide good quality solutions with substantially less computational effort.
The study is done to demonstrate the computational efficiency of the greedy heuristic over a solver, and therefore we only tabulate results for the RA11C spread here. Similar results for TN11C and RAEPC can be found in Appendix~\ref{sec:fn_compstudy_walltime_appendix}. 

\subsection{Evaluating quality of greedy solution using MRP}
\label{sec:fn_compstudy_greedysolutionquality}

Next, we evaluate the solutions from the greedy heuristic using the multiple replications procedure (MRP) described in \citet{bayraksan2006assessing}. The idea is to evaluate these candidate solutions on a larger set of independent sample paths and estimate the optimality gap using the upper bound obtained from solving a linear programming relaxation of SAA on this larger set. This is then repeated $n_g$ times to form an absolute gap estimate $\overline{G}_n(n_g)$ and the corresponding one-sided approximate (1-$\alpha$)-level confidence interval $[0, \overline{G}_n(n_g) + \epsilon_g]$. We refer readers to \cite{bayraksan2006assessing} for full details.

\begin{figure}
    \centering
    \begin{tikzpicture}
        \begin{axis}[axis lines=left, xlabel=$n$, ylabel={Gap estimate ($\overline{G}_n$)}, width=0.6\textwidth, ymajorgrids=true, xmajorgrids=true, grid style=dotted, legend style={at={(0.98, 0.98)}, anchor=north east, legend cell align=left}, xtick=data, xmin=0, ytick={0, 0.02, 0.04, 0.06, 0.08, 0.1}, ymin=0, ymax=0.11, font=\footnotesize]
            \addplot[color=graphred,mark=triangle,line width=1pt,mark size=2.5pt] table[x=reps, y=TN11C_k100_t4_r5_mean, col sep=comma] {filecontents/gapestimate.csv};
            \addlegendentry{TN11C: $k=100, t_0=4, r=0.05$};
            \addplot[color=graphblue,mark=o,line width=1pt,mark size=2.5pt] table[x=reps, y=RA1PC_k50_t3_r5_mean, col sep=comma] {filecontents/gapestimate.csv};
            \addlegendentry{RA11C: $k=50, t_0=3, r=0.05$};
        \end{axis}
    \end{tikzpicture}
    \caption{Plot of $\overline{G}_n$ vs.\ $n$ (the number of samples used to get solution using greedy heuristic)}
    \label{fig:gapestimate}
\end{figure}
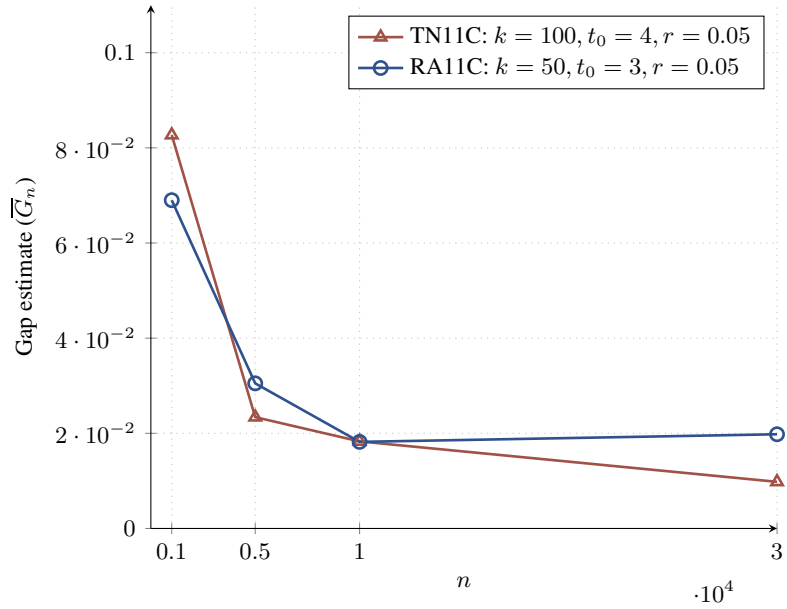

In our study we use 50,000 virus sample paths to evaluate the candidate solutions in each of $n_g=20$ replications, and $\alpha=0.05$ to form the one-sided confidence interval (CI). Due to the computational effort required for this procedure, we perform experiments only for a limited set of parameters. The results for the TN11C and the RA1PC models appear in Table~\ref{tab:evaluate_greedy_solution_gap_estimate_TN11C_k100_t4_r5} and Table~\ref{tab:evaluate_greedy_solution_gap_estimate_RA1PC_k50_t3_r5}, respectively. We plot the gap estimate $\overline{G}_n(n_g)$ for the TN11C and the RA11C spread models in Figure~\ref{fig:gapestimate}. As expected, the gap generally  decreases with an increase in the number of samples used to generate the greedy the solution, i.e., the greedy solution is expected to improve with larger sample sizes. We also note that in these tests, the greedy heuristic performs quite well, generating solutions within about 2\% of optimal for when 5000 samples or more are used.


\begin{table}[ptbh]
    \centering
    \caption{Assessing the quality of greedy solutions using MRP for 6-core EU email network under TN11C with $k=100$, $t_0=4$ and $r=0.05$. For varying values of sample size $n$, we tabulate the gap estimate $\overline{G}_n$ of the greedy solution and the width $\epsilon_g$ of the corresponding one-sided CI for $n_g=20$ and $\alpha=0.05$.}
    \label{tab:evaluate_greedy_solution_gap_estimate_TN11C_k100_t4_r5}
    \begin{tabular}{@{}ccc@{}}
        \toprule
        $\bm{n}$ & $\bm{\overline{G}_n(n_g)}$ & $\bm{\epsilon_g}$ \\ \midrule
        1000 & 0.0827 & 0.0008 \\
        5000 & 0.0234 & 0.0006 \\
        10000 & 0.0183 & 0.0005 \\
        30000 & 0.0098 & 0.0004 \\ 
        \bottomrule
    \end{tabular}
\end{table}

\begin{table}[ptbh]
    \centering
    \caption{Assessing the quality of greedy solutions using MRP for 6-core EU email network under RA11C with $k=50$, $t_0=3$ and $r=0.05$. For varying values of sample size $n$, we tabulate the gap estimate $\overline{G}_n$ of the greedy solution and the width $\epsilon_g$ of the corresponding one-sided CI for $n_g=20$ and $\alpha=0.05$.}
    \label{tab:evaluate_greedy_solution_gap_estimate_RA1PC_k50_t3_r5}
    \begin{tabular}{@{}ccc@{}}
        \toprule
        $\bm{n}$ & $\bm{\overline{G}_n(n_g)}$ & $\bm{\epsilon_g}$ \\ \midrule
        1000 & 0.0690 & 0.0006 \\ 
        5000 & 0.0305 & 0.0006 \\ 
        10000 & 0.0182 & 0.0005 \\ 
        30000 & 0.0198 & 0.0006 \\ 
        \bottomrule
    \end{tabular}
\end{table}

\subsection{Effect of false-negative rate on the probability of detection}
\label{sec:fn_compstudy_affect_of_r}

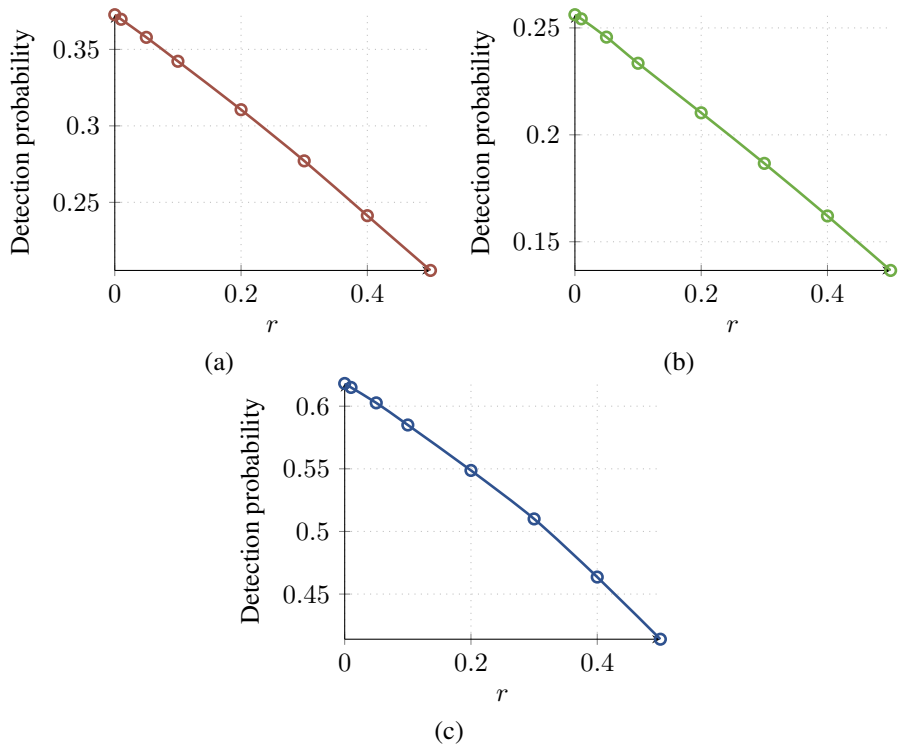
\begin{figure}[ptbh]
    \centering
    $
        \begin{array}{cc}
            \begin{tikzpicture}
                \begin{axis}[axis lines=left, xlabel=$r$, ylabel=Detection probability, width=0.35\textwidth, ymajorgrids=true, xmajorgrids=true, grid style=dotted, tick=data]
                    \addplot[smooth, color=graphred,mark=o,line width=1pt,mark size=2pt] table[x=r, y=objective, col sep=comma]{filecontents/TN11C_heuristic_point_estimate.csv};
                \end{axis}
            \end{tikzpicture}
            & 
            \begin{tikzpicture}
                \begin{axis}[axis lines=left, xlabel=$r$, ylabel=Detection probability, width=0.35\textwidth, ymajorgrids=true, xmajorgrids=true, grid style=dotted, tick=data]
                    \addplot[smooth, color=graphgreen,mark=o,line width=1pt,mark size=2pt] table[x=r, y=objective, col sep=comma]{filecontents/RA1PC_heuristic_point_estimate.csv};
                \end{axis}
            \end{tikzpicture}
            \\
            \text{(a)} & \text{(b)} \\
            \multicolumn{2}{c}
            {
            \begin{tikzpicture}
                \begin{axis}[axis lines=left, xlabel=$r$, ylabel=Detection probability, width=0.35\textwidth, ymajorgrids=true, xmajorgrids=true, grid style=dotted, tick=data]
                    \addplot[smooth, color=graphblue,mark=o,line width=1pt,mark size=2pt] table[x=r, y=objective, col sep=comma]{filecontents/RAEPC_heuristic_point_estimate.csv};
                \end{axis}
            \end{tikzpicture} 
            } \\
            \multicolumn{2}{c}{\text{(c)}}
        \end{array}
    $
    \caption{Effect of increasing false-negative rate on probability of detection for $k=50$ and $n=50000$: (a) TN11C spread, (b) RA1PC spread with $p=0.5$ and (c) RAEPC spread with $p=0.5$}
    \label{fig:heuristic_point_estimate}
\end{figure}

Having established the computational efficacy of the greedy heuristic and the high quality of solutions produced by it, we now use greedy solutions to study the effect of varying the false negative probability $r$ (from 0 to 0.5) on the probability of detection. We use the greedy heuristic to get an ``approximately optimal'' solution on 50,000 samples for $k=50$ (about 1\% of the total number of nodes), and then evaluate this solution on an independent set of 2 million samples to get a point estimate of the 
approximate optimal objective value, for different spread dynamics and varying values of $r$. We plot the results for the TN11C, the RA1PC, and the RAEPC models in Figure~\ref{fig:heuristic_point_estimate}. The probability of detection seems to drop approximately linearly in $r$ for all the three cases. For TN11C the probability drops by 0.17 when using a 50\% reliable detector, a 45\% decrease from the overall probability of detection with a perfectly reliable detector. For RA1PC the probability drops by 0.12, a 47\% decrease, and for RAEPC the probability drops by 0.20, a 33\% decrease. For a 95\% reliable detector (i.e., $r=0.05$) the drop is around 4\% for TN11C and RA1PC, and 3\% for RAEPC.

\subsection{Ignoring Detector Fallibility}
\label{sec:fn_costmodelingfn}
In our final computational study, we examine the effect on solution quality when the decision maker ignores the fallibility of detectors. In particular, we solve a model that assumes that the detectors are completely reliable, and compare it to the solution of the ``true'' model, with fallible detectors.  

The following procedure is used to generate the data in Table~\ref{tab:semihamming_TN11C_t4}, Table~\ref{tab:semihamming_RA11C_t3}, and Table~\ref{tab:semihamming_RAEPC_t1_p50}:
\begin{enumerate}[leftmargin=0.25in]
    \item Generate 50,000 sample paths of virus spread. \label{enum:fn_costmodelingfn_1}
    \item Generate a VDM (i.e., coin flips to model false negatives) corresponding to each of the 50,000 sample paths in Step 1. The 50,000 sample paths generated in Step~\ref{enum:fn_costmodelingfn_1} along with corresponding VDMs serve as 50,000 samples for the ``true model'' i.e., the model with false-negative results. \label{enum:fn_costmodelingfn_2}
    \item Generate 5 million sample paths, independently from the Step~\ref{enum:fn_costmodelingfn_1} paths, for virus spread, along with the corresponding VDMs. \label{enum:fn_costmodelingfn_3}
    \item Use the greedy heuristic on sample paths generated in Step 1, to select nodes at which to place detectors (ignoring false negatives). \label{enum:fn_costmodelingfn_4}
    \item Use the greedy heuristic on the 50,000 samples in Step~\ref{enum:fn_costmodelingfn_1} and the VDMs in Step~\ref{enum:fn_costmodelingfn_2}, corresponding to the true model, to select nodes at which to place detectors. \label{enum:fn_costmodelingfn_5}
    \item Evaluate the solutions from Step~\ref{enum:fn_costmodelingfn_4} and Step~\ref{enum:fn_costmodelingfn_5} on the samples generated in Step~\ref{enum:fn_costmodelingfn_3}. The corresponding objective function values are referred to as FN1 and FN2, respectively, in Table~\ref{tab:semihamming_TN11C_t4} through Table~\ref{tab:semihamming_RAEPC_t1_p50}. Specifically, we calculate (a) the difference in the probability of detection and (b) half of the Hamming distance (the semi-hamming distance), a measure of how different the detector node sets are. \label{enum:fn_costmodelingfn_6}
\end{enumerate}

\begin{table}[ptbh]
    \centering
    \caption{Semi-hamming distance for TN11C spread with $t_0=4$ and varying false-negative rates}
    \label{tab:semihamming_TN11C_t4}
    \begin{tabular}{@{}cccccc@{}}
    \toprule
    $\bm{k}$ & $\bm{r}$ & \textbf{\begin{tabular}[c]{@{}c@{}}Objective\\ ($\bm{r=0}$)\\ FN1\end{tabular}} & \textbf{\begin{tabular}[c]{@{}c@{}}Objective\\ (true model)\\ FN2\end{tabular}} & \textbf{\begin{tabular}[c]{@{}c@{}}difference \\ (FN1-FN2)\end{tabular}} & \textbf{\begin{tabular}[c]{@{}c@{}}half of \\ hamming \\ distance\end{tabular}} \\ \midrule
    \multirow{4}{*}{100} & 0.05 & 0.5405 & 0.5407 & -0.0002 & 3 \\
     & 0.1 & 0.5212 & 0.5217 & -0.0005 & 3 \\
     & 0.25 & 0.4564 & 0.4549 & 0.0015 & 8 \\
     & 0.3 & 0.4312 & 0.4323 & -0.0011 & 7 \\ \cmidrule(l){2-6}
    \multirow{4}{*}{200} & 0.05 & 0.7463 & 0.7465 & -0.0002 & 6 \\
     & 0.1 & 0.7250 & 0.7247 & 0.0003 & 8 \\
     & 0.25 & 0.6480 & 0.6488 & -0.0008 & 12 \\
     & 0.3 & 0.6183 & 0.6209 & -0.0026 & 20 \\ \cmidrule(l){2-6}
    \multirow{4}{*}{250} & 0.05 & 0.8079 & 0.8083 & 0.0311 & 12 \\
     & 0.1 & 0.7865 & 0.7875 & -0.0014 & 16 \\
     & 0.25 & 0.7110 & 0.7116 & 0.0045 & 21 \\
     & 0.3 & 0.6831 & 0.6852 & -0.0006 & 19 \\ \cmidrule(l){2-6} 
    \end{tabular}
\end{table}

\begin{table}[ptbh]
    \centering
    \caption{Semi-hamming distance for RA11C spread with $t_0=3$ and varying false-negative rates}
    \label{tab:semihamming_RA11C_t3}
    \begin{tabular}{@{}cccccc@{}}
    \toprule
    $\bm{k}$ & $\bm{r}$ & \textbf{\begin{tabular}[c]{@{}c@{}}Objective\\ ($\bm{r=0}$)\\ FN1\end{tabular}} & \textbf{\begin{tabular}[c]{@{}c@{}}Objective\\ (true model)\\ FN2\end{tabular}} & \textbf{\begin{tabular}[c]{@{}c@{}}difference \\ (FN1-FN2)\end{tabular}} & \textbf{\begin{tabular}[c]{@{}c@{}}half of \\ hamming \\ distance\end{tabular}} \\ \midrule
    \multirow{4}{*}{50} & 0.05 & 0.5198 & 0.5187 & 0.0011 & 1 \\
     & 0.1 & 0.4992 & 0.4991 & 0.0001 & 2 \\
     & 0.25 & 0.4363 & 0.4356 & 0.0007 & 4 \\
     & 0.3 & 0.4137 & 0.4136 & 0.0001 & 4 \\ \cmidrule(l){2-6}
    \multirow{4}{*}{100} & 0.05 & 0.7228 & 0.7221 & 0.0007 & 4 \\
     & 0.1 & 0.7023 & 0.7020 & 0.0003 & 4 \\
     & 0.25 & 0.6277 & 0.6284 & -0.0007 & 10 \\
     & 0.3 & 0.6006 & 0.6024 & -0.0018 & 10 \\ \cmidrule(l){2-6}
    \multirow{4}{*}{200} & 0.05 & 0.8898 & 0.8893 & 0.0005 & 15 \\
     & 0.1 & 0.8729 & 0.8737 & -0.0008 & 18 \\
     & 0.25 & 0.8100 & 0.8130 & -0.0030 & 28 \\
     & 0.3 & 0.7850 & 0.7883 & -0.0033 & 24 \\ \cmidrule(l){2-6} 
    \end{tabular}
\end{table}

\begin{table}[ptbh]
    \centering
    \caption{Semi-hamming distance for RAEPC spread with $t_0=1$, $p=0.5$ and varying false-negative rates}
    \label{tab:semihamming_RAEPC_t1_p50}
    \begin{tabular}{@{}cccccc@{}}
    \toprule
    $\bm{k}$ & $\bm{r}$ & \textbf{\begin{tabular}[c]{@{}c@{}}Objective\\ ($\bm{r=0}$)\\ FN1\end{tabular}} & \textbf{\begin{tabular}[c]{@{}c@{}}Objective\\ (true model)\\ FN2\end{tabular}} & \textbf{\begin{tabular}[c]{@{}c@{}}difference \\ (FN1-FN2)\end{tabular}} & \textbf{\begin{tabular}[c]{@{}c@{}}half of \\ hamming \\ distance\end{tabular}} \\ \midrule
    \multirow{4}{*}{50} & 0.05 & 0.6022 & 0.6023 & -0.0001 & 2 \\
     & 0.10 & 0.5865 & 0.5865 & 0.0000 & 2 \\
     & 0.25 & 0.5297 & 0.5306 & -0.0007 & 2 \\ 
     & 0.30 & 0.5074 & 0.5085 & -0.0011 & 7 \\ \cmidrule(l){2-6}
    \multirow{4}{*}{100} & 0.05 & 0.7564 & 0.7555 & 0.0009 & 7 \\ 
     & 0.10 & 0.7379 & 0.7400 & -0.0021 & 5 \\
     & 0.25 & 0.6841 & 0.6839 & 0.0002 & 7 \\
     & 0.30 & 0.6623 & 0.6636 & -0.0013 & 12 \\ \cmidrule(l){2-6}
    \end{tabular}
\end{table}

The semi-hamming distance between two sets of detector nodes indicates the minimum number of node changes that are required to make the two sets the same. From Table~\ref{tab:semihamming_TN11C_t4} through Table~\ref{tab:semihamming_RAEPC_t1_p50} we make the following observations:
\begin{enumerate}[leftmargin=0.25in]
    \item There seems to be a negligible absolute difference in the probability of detection.
    \item For a fixed $k$, the semi-Hamming distance generally increases but is not always monotonically increasing in $r$.
    \item For $r=0.3$, the semi-Hamming distance is as high as 10\% of the number of detectors.
\end{enumerate}

The results above induce at least two questions regarding modeling of false negatives. First, why does a model that ignores false negatives select a detector set which is close to the ``correct'' set from the true model, at least when using a greedy algorithm? 
Second, why do differences in up to two dozen nodes induce very small changes in objective values? The likely answer to the second question is simply that as long as a ``core'' set of nodes is selected correctly, changing the remaining 5\%--10\% of nodes in the detector set does not matter because of small marginal detection gains as the last nodes are added. As for the first question, it is harder to guess the root cause of this effect.


From the results above, one might infer that properly modeling false negatives is not particularly important in our models. However, we should distinguish two different purposes for building such models. First, these models are useful in \emph{design} decisions, i.e., deciding where detectors should be placed. Second, once locations are chosen, the models are used for \emph{performance} modeling, e.g., evaluating the probability of successful detection for a given design. Our somewhat limited results indicate that modeling false negatives may not be crucial in the design stage. However, in our performance evaluation in Section~\ref{sec:fn_compstudy_affect_of_r} there was a  significant decrease in detection probability with increases in the false-negative rate. Hence, it is important to incorporate false negatives when making a decision, for example, on how many detectors are needed to ensure a certain probability of detection.


\section{Optimal Placement in a Wheel Network}
\label{sec:fn_wheelnetworkplacement}

\begin{figure}[ptbh]
    \centering
    $
        \begin{array}{cc}
            \scalebox{0.6}
            {
                \begin{tikzpicture}
                    \draw[color=graphgreen, line width=1pt](0,0) circle (1.5in);
                    \draw[color=graphblue, fill=graphblue!80](0,0) circle (0.1in);
                    \draw[color=graphblue, fill=graphblue!80](1.5in,0) circle (0.1in);
                    \draw[color=graphblue, fill=graphblue!80](-1.5in,0) circle (0.1in);
                    \draw[color=graphblue, fill=graphblue!80](0, 1.5in) circle (0.1in);
                    \draw[color=graphblue, fill=graphblue!80](0, -1.5in) circle (0.1in);
                    \draw[color=graphblue, fill=graphblue!80](1.06in,1.06in) circle (0.1in);
                    \draw[color=graphblue, fill=graphblue!80](-1.06in,1.06in) circle (0.1in);
                    \draw[color=graphblue, fill=graphblue!80](1.06in,-1.06in) circle (0.1in);
                    \draw[color=graphblue, fill=graphblue!80](-1.06in,-1.06in) circle (0.1in);
                    \draw[color=graphgreen!80, line width=1pt] (0.1in,0) -- (1.4in,0);
                    \draw[color=graphgreen!80, line width=1pt] (-0.1in,0) -- (-1.4in,0);
                    \draw[color=graphgreen!80, line width=1pt] (0,0.1in) -- (0, 1.4in);
                    \draw[color=graphgreen!80, line width=1pt] (0,-0.1in) -- (0, -1.4in);
                    \draw[color=graphgreen!80, line width=1pt] (0.07in,0.07in) -- (0.99in,0.99in);
                    \draw[color=graphgreen!80, line width=1pt] (-0.07in,0.07in) -- (-0.99in,0.99in);
                    \draw[color=graphgreen!80, line width=1pt] (0.07in,-0.07in) -- (0.99in,-0.99in);
                    \draw[color=graphgreen!80, line width=1pt] (-0.07in,-0.07in) -- (-0.99in,-0.99in);
                    \node (detector1) [detector, yshift=1.5in] {D1};
                    \node (detector2) [detector, xshift=1.06in, yshift=1.06in] {D2};
                \end{tikzpicture}
            }
            & 
            \scalebox{0.6}
            {
                \begin{tikzpicture}
                    \draw[color=graphgreen, line width=1pt](0,0) circle (1.5in);
                    \draw[color=graphblue, fill=graphblue!80](0,0) circle (0.1in);
                    \draw[color=graphblue, fill=graphblue!80](1.5in,0) circle (0.1in);
                    \draw[color=graphblue, fill=graphblue!80](-1.5in,0) circle (0.1in);
                    \draw[color=graphblue, fill=graphblue!80](0, 1.5in) circle (0.1in);
                    \draw[color=graphblue, fill=graphblue!80](0, -1.5in) circle (0.1in);
                    \draw[color=graphblue, fill=graphblue!80](1.06in,1.06in) circle (0.1in);
                    \draw[color=graphblue, fill=graphblue!80](-1.06in,1.06in) circle (0.1in);
                    \draw[color=graphblue, fill=graphblue!80](1.06in,-1.06in) circle (0.1in);
                    \draw[color=graphblue, fill=graphblue!80](-1.06in,-1.06in) circle (0.1in);
                    \draw[color=graphgreen!80, line width=1pt] (0.1in,0) -- (1.4in,0);
                    \draw[color=graphgreen!80, line width=1pt] (-0.1in,0) -- (-1.4in,0);
                    \draw[color=graphgreen!80, line width=1pt] (0,0.1in) -- (0, 1.4in);
                    \draw[color=graphgreen!80, line width=1pt] (0,-0.1in) -- (0, -1.4in);
                    \draw[color=graphgreen!80, line width=1pt] (0.07in,0.07in) -- (0.99in,0.99in);
                    \draw[color=graphgreen!80, line width=1pt] (-0.07in,0.07in) -- (-0.99in,0.99in);
                    \draw[color=graphgreen!80, line width=1pt] (0.07in,-0.07in) -- (0.99in,-0.99in);
                    \draw[color=graphgreen!80, line width=1pt] (-0.07in,-0.07in) -- (-0.99in,-0.99in);
                    \node (detector1) [detector, yshift=1.5in] {D1};
                    \node (detector2) [detector, yshift=-1.5in] {D3};
                \end{tikzpicture}
            }
            \\
            \text{(a)} & \text{(b)} \\
            \scalebox{0.6}
            {
                \begin{tikzpicture}
                    \draw[color=graphgreen, line width=1pt](0,0) circle (1.5in);
                    \draw[color=graphblue, fill=graphblue!80](0,0) circle (0.1in);
                    \draw[color=graphblue, fill=graphblue!80](1.5in,0) circle (0.1in);
                    \draw[color=graphblue, fill=graphblue!80](-1.5in,0) circle (0.1in);
                    \draw[color=graphblue, fill=graphblue!80](0, 1.5in) circle (0.1in);
                    \draw[color=graphblue, fill=graphblue!80](0, -1.5in) circle (0.1in);
                    \draw[color=graphblue, fill=graphblue!80](1.06in,1.06in) circle (0.1in);
                    \draw[color=graphblue, fill=graphblue!80](-1.06in,1.06in) circle (0.1in);
                    \draw[color=graphblue, fill=graphblue!80](1.06in,-1.06in) circle (0.1in);
                    \draw[color=graphblue, fill=graphblue!80](-1.06in,-1.06in) circle (0.1in);
                    \draw[color=graphgreen!80, line width=1pt] (0.1in,0) -- (1.4in,0);
                    \draw[color=graphgreen!80, line width=1pt] (-0.1in,0) -- (-1.4in,0);
                    \draw[color=graphgreen!80, line width=1pt] (0,0.1in) -- (0, 1.4in);
                    \draw[color=graphgreen!80, line width=1pt] (0,-0.1in) -- (0, -1.4in);
                    \draw[color=graphgreen!80, line width=1pt] (0.07in,0.07in) -- (0.99in,0.99in);
                    \draw[color=graphgreen!80, line width=1pt] (-0.07in,0.07in) -- (-0.99in,0.99in);
                    \draw[color=graphgreen!80, line width=1pt] (0.07in,-0.07in) -- (0.99in,-0.99in);
                    \draw[color=graphgreen!80, line width=1pt] (-0.07in,-0.07in) -- (-0.99in,-0.99in);
                    \node (detector1) [detector, yshift=1.5in] {D1};
                    \node (detector2) [detector] {D5};
                \end{tikzpicture}
            } 
            & 
            \scalebox{0.6}
            {
                \begin{tikzpicture}
                    \draw[color=graphgreen, line width=1pt](0,0) circle (1.5in);
                    \draw[color=graphblue, fill=graphblue!80](0,0) circle (0.1in);
                    \draw[color=graphblue, fill=graphblue!80](1.5in,0) circle (0.1in);
                    \draw[color=graphblue, fill=graphblue!80](-1.5in,0) circle (0.1in);
                    \draw[color=graphblue, fill=graphblue!80](0, 1.5in) circle (0.1in);
                    \draw[color=graphblue, fill=graphblue!80](0, -1.5in) circle (0.1in);
                    \draw[color=graphblue, fill=graphblue!80](1.06in,1.06in) circle (0.1in);
                    \draw[color=graphblue, fill=graphblue!80](-1.06in,1.06in) circle (0.1in);
                    \draw[color=graphblue, fill=graphblue!80](1.06in,-1.06in) circle (0.1in);
                    \draw[color=graphblue, fill=graphblue!80](-1.06in,-1.06in) circle (0.1in);
                    \draw[color=graphgreen!80, line width=1pt] (0.1in,0) -- (1.4in,0);
                    \draw[color=graphgreen!80, line width=1pt] (-0.1in,0) -- (-1.4in,0);
                    \draw[color=graphgreen!80, line width=1pt] (0,0.1in) -- (0, 1.4in);
                    \draw[color=graphgreen!80, line width=1pt] (0,-0.1in) -- (0, -1.4in);
                    \draw[color=graphgreen!80, line width=1pt] (0.07in,0.07in) -- (0.99in,0.99in);
                    \draw[color=graphgreen!80, line width=1pt] (-0.07in,0.07in) -- (-0.99in,0.99in);
                    \draw[color=graphgreen!80, line width=1pt] (0.07in,-0.07in) -- (0.99in,-0.99in);
                    \draw[color=graphgreen!80, line width=1pt] (-0.07in,-0.07in) -- (-0.99in,-0.99in);
                    \node (detector1) [detector, minimum width=0.5cm, minimum height=1cm] {D5};
                    \node (detector2) [detector] {D5};
                \end{tikzpicture}
            } 
            \\
            \text{(c)} & \text{(d)} \\
        \end{array}
    $
    \caption{Different detector placements: (a) Detectors D1 and D2 on adjacent-edge nodes, (b) detectors D1 and D3 on diametrically-opposite edge nodes, (c) detector D1 on an edge node and detector D5 on the center node and (d) both detectors on the center node}
    \label{fig:detector_placements}
\end{figure}
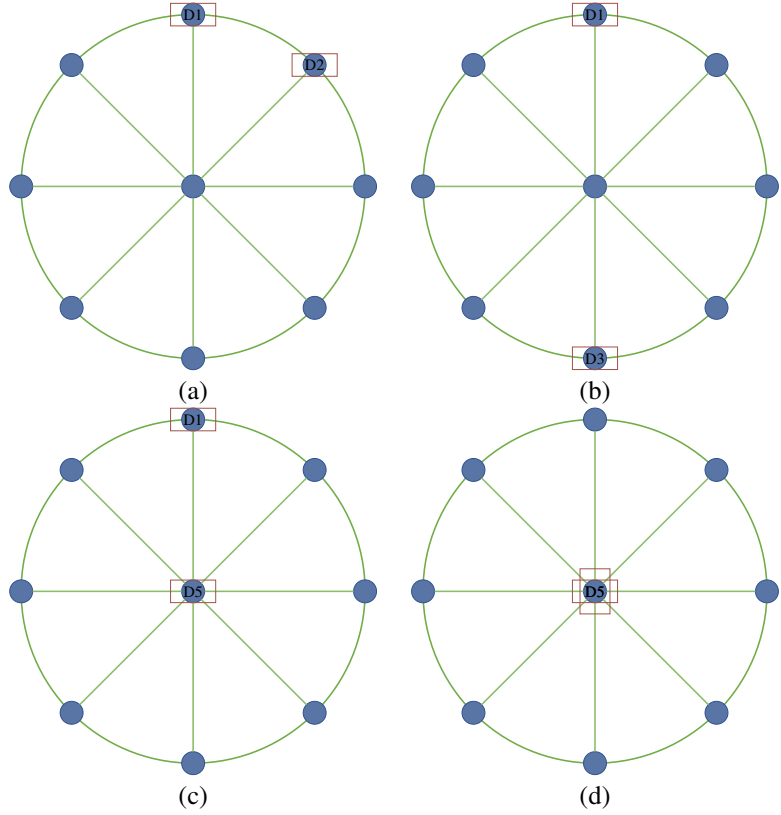
In this section we provide some analytical insight, via stylized examples, on optimal placement of detectors. Specifically, we develop analytical expressions for detector placement strategies around high-degree nodes, which are generally chosen for influence purposes in social networks. We consider the wheel network, four versions of which are depicted in Figure~\ref{fig:detector_placements}. The wheel network has one central ``high-degree'' node that connects to many degree-3 edge nodes. We derive analytical expressions for the probability of detecting a virus under the following assumptions:
\begin{itemize}[leftmargin=0.25in]
    \item There is only one time step, i.e., $t=0\rightarrow t=1$ (two ``days'' of detection).
    \item There are two detectors available. 
    \item The spread model is TN11C.
    \item The number of nodes, $v$, is greater than or equal to 5 (one center node and at least 4 edge nodes).
    \item The false negative probability satisfies $0\leq r<1$ for all detectors.
\end{itemize}

In Table~\ref{tab:detectionprobability} we give the probability of detecting the virus for the four different configurations shown in Figure~\ref{fig:detector_placements}. Note that one of the configurations places both detectors on the central node. There is a high likelihood of the virus visiting the central node in two detection steps, and one might wish to place two detectors at the same node if the detectors are unreliable. (We assume that they produce false negatives independently.) We also provide the limiting probability of detecting the virus as $v$ increases, when the false-negative rate $r$ is held constant. The derivations for these expressions in the table and proofs for our corresponding claims below appear in Appendix~\ref{sec:fn_wheelnetworkplacement_appendix}.

\begin{table}[ptbh]
    \centering
    \caption{Probability of detecting the virus for different detector placement configurations shown in Figure~\ref{fig:detector_placements}}
    \label{tab:detectionprobability}
    \begin{tabular}{lcc}
        \toprule
        \textbf{Detector placement} & \textbf{Detection probability} & \textbf{\begin{tabular}[c]{@{}c@{}}Detection probability for \\ fixed $r$, $v\rightarrow\infty$\end{tabular}} \\ \midrule
        \multirow{2}{*}{Adjacent-edge nodes} & \multirow{2}{*}{$\frac{2(1-r)}{3v(v-1)} \left[4v-1+rv-r \right]$} & \multirow{2}{*}{$0$} \\
         &  &  \\
        \multirow{2}{*}{Diametrically-opposite nodes} & \multirow{2}{*}{$\frac{2(1-r)}{3v(v-1)} \left(5v-2 \right)$} & \multirow{2}{*}{$0$} \\
         &  &  \\
        \multirow{2}{*}{Center-edge nodes} & \multirow{2}{*}{$\frac{(1-r)}{3v(v-1)}\left[v^2+5v+rv+2r-6  \right]$} & \multirow{2}{*}{$\frac{1-r}{3}$} \\
         &  &  \\
        \multirow{2}{*}{Center-center nodes} & \multirow{2}{*}{$\frac{(1-r)(1+r)(v+2)}{3v}$} & \multirow{2}{*}{$\frac{1-r^2}{3}$} \\
         &  &  \\
        \bottomrule
    \end{tabular}
\end{table}

We have the following results from the limiting values ($v\rightarrow\infty$) in Table~\ref{tab:detectionprobability}:
\begin{enumerate}[leftmargin=0.25in]
    \item For any $r >0$, placing two detectors at the central node is optimal, for sufficiently large $v$.
    \item When detectors are completely reliable, i.e., $r=0$, at least one detector is placed at the central node in an optimal solution.
\end{enumerate}
For smaller values of $v$, we observe the following:
\begin{enumerate}[leftmargin=0.25in]
    \item Placing detectors on diametrically-opposite edge nodes is always better than placing them on adjacent-edge nodes.
    \item For $v\geq 6$, the center-edge placement is better than placing the detectors on diametrically-opposite edge nodes. For $v=5$,  if $r>\frac{2}{7}$, then the center-edge placement is better, and if $r<\frac{2}{7}$, then the diametrically-opposite placement is better.
    \item For a given $r$, let $\left\lceil v_\text{min} \right\rceil$ be the minimum number of nodes needed for the center-center placement to be better than the center-edge placement of two detectors. Table~\ref{tab:comparisoncentercenter_centeredgenprime} provides $\left\lceil v_\text{min} \right\rceil$ for a few reasonable values of $r$.
     We see that for higher false-negative rates it might be better to place two detectors at the central node in a cluster. 
    \item For $r \geq \frac{16}{21} \approx 0.76$, the center-center placement is at least as good as the center-edge placement, because $\left\lceil v_\text{min} \right\rceil = 5$ at $r= \frac{16}{21}$. 
\end{enumerate}

\begin{table}[ptbh]
    \centering
    \caption{Minimum number of nodes required for center-center placement of detectors to be better than center-edge placement of detectors, for some standard values of $r$}
    \label{tab:comparisoncentercenter_centeredgenprime}
    \begin{tabular}{cc}
        \toprule
        $\bm{r}$ & $\bm{\left\lceil v_\text{min} \right\rceil}$ \\ \midrule
        0.01 & 400 \\ 
        0.05 & 80 \\ 
        0.10 & 40 \\ 
        \bottomrule
    \end{tabular}
\end{table}
Analytical results for larger values of $t$ can be derived, but the expressions become intractable quickly, even for moderate values of $t$.

\section{Final Remarks}
\label{sec:conclusion}


A natural extension to this research is to allow for detectors with false-positive results. It is clear that this extension requires a substantially modified optimization model. In the false-positive case the network manager has to make two decisions: a design decision and an operational decision. This is unlike the false-negative case where any false negative does not trigger an ``alarm,'' and the virus just goes undetected. The design problem is to decide where to place the honeypots, at time 0. The operational problem is what action to take when a honeypot signals a virus detection. For example, the network manager can choose to ignore the alarm, in which case a possible virus detection is missed, but a false positive is avoided.
Hence, in a model with false positives, we seek an optimal operational policy that incorporates both the costs of taking action when there is false positive, and a failure to detect the virus.  
It would be interesting to see if the induced Markov decision problem can be formulated as an integer program. It is not apparent if the submodularity property would still hold, however it might be possible to develop a simulation model with a ``similar'' greedy heuristic and gain some insights into the problem through computational results. One could also explore how to combine all these models and develop a generalized model which accounts for the possibility of detectors being unreliable and producing both false-positive results and false-negative results.

Another possible extension to our work is to model the choices that malicious players (adversaries) make in order to inflict maximum damage to the system. The adversary attempts to maximize the expected number of infected nodes by the time the virus is detected by introducing the virus on a single node, and the network manager installs honeypots to minimize this value. It would be interesting to explore this min-max model extension.

\section*{Acknowledgements}
{
    The authors would like to acknowledge the feedback of Dave Morton and Grani Hanasusanto which helped in refining the exposition of this research.
}

\appendix
\appendixpage
\section{EU email network}
\label{sec:euemailnetwork_appendix}
\begin{defn}[\emph{$c$-core}]\label{defn:ccore}
    Given a graph $G = (V, E)$, $G_S = (V_S, E_S)$ is a subgraph of $G$ that has $V_S \subseteq V$ and $E_S = \{(i,j)\in E: i\in V_S, j\in V_S \}$. A subgraph $G_S = (V_S, E_S)$ induced by the set $V_S \subset V$ is a $c$-core if and only if $G_S$ is a maximum subgraph for which the degree of every vertex $v\in V_S$ in $G_S$ is greater than or equal to $c$. A maximum induced subgraph has as many vertices as possible.
\end{defn}

\begin{table}[ptbh]
    \centering
    \caption{Some network characteristics of the largest connected component of several $c$-core decompositions of the EU email communication network}
    \label{tab:eunetwork_characteristics_appendix}
    \begin{tabular}{cllcc}
        \toprule
        $\bm{c}$ & \multicolumn{1}{c}{\textbf{\#nodes}} & \multicolumn{1}{c}{\textbf{\#edges}} & \textbf{\begin{tabular}[c]{@{}c@{}}max.\\ degree\end{tabular}} & \textbf{\begin{tabular}[c]{@{}c@{}}avg.\\ degree\end{tabular}} \\ \midrule
        1 & 224832 & 339925 & 7636 & 3.02 \\ 
        2 & 36167 & 151260 & 1099 & 8.36 \\ 
        4 & 9914 & 92897 & 703 & 18.7 \\ 
        6 & 5400 & 73315 & 552 & 27.2 \\ 
        8 & 3654 & 62135 & 476 & 34.0 \\ 
        \bottomrule
    \end{tabular}
\end{table}

\section{Computational Studies}
\label{sec:fn_computationalstudies_appendix}

\subsection{Choosing a time threshold}
\label{sec:fn_choosingtime_appendix}
We generate Monte Carlo samples to determine the average time it takes a fixed percentage of nodes to become infected. This analysis provided guidance in selecting the value of $t_0$ in our computational studies. We selected infection of 0.1\% of the nodes as the threshold, but a decision maker may of course adjust this threshold. In Table~\ref{tab:estimating_t0_appendix} we present the results of such simulations for different levels of infections. These results were obtained by simulating the virus spread using TN11C, RA11C and RAE1C dynamics, terminating the simulation when the given percentage of infections was achieved, and recording the time. The average time is computed over 10,000 independent replications. Since the initial infection happens at $t=0$, if a virus is detected at, say $t=3$, then it means the virus has been detected within 4 time slots. 

\begin{table}[ptbh]
    \centering
    \caption{Simulation results estimating a time threshold $t_0$ for TN11C, RA11C and RAE1C, through 10000 simulation repetitions of virus spread on the 6-core EU email network}
    \label{tab:estimating_t0_appendix}
    \begin{tabular}{cccc}
        \toprule
        \textbf{\begin{tabular}[c]{@{}c@{}}Percent \\ infection\end{tabular}} & \textbf{\begin{tabular}[c]{@{}c@{}}Average \\ time for \\ TN11C\end{tabular}} & \textbf{\begin{tabular}[c]{@{}c@{}}Average \\ time for \\ RA11C\end{tabular}} & \textbf{\begin{tabular}[c]{@{}c@{}}Average \\ time for \\ RAE1C\end{tabular}} \\ \midrule
        0.1 & 4.11 & 3.00 & 1 \\ 
        0.5 & 27.3 & 5.39 & 1.78 \\ 
        1.0 & 56.3 & 6.62 & 1.88 \\ 
        5.0 & 309 & 9.27 & 2.02 \\ 
        \bottomrule
    \end{tabular}
    
\end{table}

In Table~\ref{tab:estimating_t0_appendix} we observe that the average time for TN11C increases approximately linearly with percent infection. In the TN11C model, the infection percentage cannot be superlinear with respect to time, because $|P_l| \leq t_0 +1, \,\forall l$. 

\subsection{Comparing MIP and greedy solutions}
\label{sec:fn_compstudy_walltime_appendix}
In Section~\ref{sec:fn_compstudy_comparemipandgreedy} we compare the MIP and greedy solutions for the RA11C model. In Table~\ref{tab:tn11c_compare_mip_and_greedy_t4_r0_appendix} and Table~\ref{tab:raepc_compare_mip_and_greedy_t1_p95_r1_appendix}, we provide comparisons for the TN11C and the RAEPC spread models, respectively. 

\begin{table}[ptbh]
    \centering
    \caption{Comparison of MIP and Greedy, under TN11C spread on 6-core EU email network with $t_0=4$ and completely reliable detectors ($r=0$): computational time and 95\% confidence interval on the difference $\mu = f(S^h_n) - f(S^\ast_n)$ using McNemar's procedure ($n^{\prime\prime}=2000000$). A positive difference indicates that the greedy solution performed better. Note that columns $\bm{\underline{z}^{\ast}_n}$, $\bm{\overline{z}^{\ast}_n}$, $\bm{z^h_n}$ and both wall times are not from the McNemar's procedure.}
    \label{tab:tn11c_compare_mip_and_greedy_t4_r0_appendix}
    \begin{tabular}{@{}cccccccc@{}}
        \toprule
        \multicolumn{2}{c}{\textbf{}} & \multicolumn{3}{|c|}{\textbf{MIP}} & \multicolumn{2}{c|}{\textbf{Greedy}} & $\bm{f(S^h_n) - f(S^\ast_n)}$ \\ \midrule
        $\bm{k}$ & $\bm{n}$ & $\bm{\underline{z}^{\ast}_n}$ & $\bm{\overline{z}^{\ast}_n}$ & \textbf{\begin{tabular}[c]{@{}c@{}}Wall time\\ (sec)\end{tabular}} & $\bm{z^h_n}$ & \textbf{\begin{tabular}[c]{@{}c@{}}Wall time\\ (sec)\end{tabular}} & \textbf{95\% CI} \\ \midrule
         & 1000 & 0.6780 & 0.6780 & 0.0583 & 0.674 & 0.1288 & {[}0.0030, 0.0037{]} \\
         & 5000 & 0.5970 & 0.5984 & 2.154 & 0.5958 & 0.1514 & {[}-0.0017, -0.0010{]} \\
         & 10000 & 0.5785 & 0.5834 & 3.132 & 0.5805 & 0.2460 & {[}0.0016, 0.0023{]} \\
         & 30000 & 0.5727 & 0.5781 & 22.22 & 0.5737 & 0.6366 & {[}-0.0006, -0.0002{]} \\
        \multirow{-5}{*}{100} & 50000 & 0.5668 & 0.5716 & 349.0 & 0.5698 & 1.008 & {[}0.0013, 0.0019{]} \\ \cmidrule(l){2-8} 
         & 1000 & 0.9170 & 0.9250 & 0.0913 & 0.9070 & 0.1774 & {[}0.0181, 0.0191{]} \\
         & 5000 & 0.8110 & 0.8188 & 14.00 & 0.8112 & 0.2434 & {[}-0.0031, -0.0023{]} \\
         & 10000 & 0.7937 & 0.7999 & 116.1 & 0.7923 & 0.4022 & {[}0.0034, 0.0039{]} \\
         & 30000 & 0.7808 & 0.7873 & 1562 & 0.7807 & 1.002 & {[}0.0020, 0.0025{]} \\
        \multirow{-5}{*}{200} & 50000 & \cellcolor[HTML]{EFEFEF}0.7542 & \cellcolor[HTML]{EFEFEF}0.7853 & \cellcolor[HTML]{EFEFEF}3600 & 0.7759 & 1.618 & \colorbox{tablegraybg}{{[}0.0198, 0.0208{]}}\\ \cmidrule(l){2-8} 
         & 1000 & 0.9790 & 0.9860 & 0.0818 & 0.9630 & 0.2110 & {[}0.0311, 0.0323{]} \\
         & 5000 & 0.8752 & 0.8832 & 65.80 & 0.8716 & 0.2964 & {[}-0.0014, -0.0008{]} \\
         & 10000 & 0.8552 & 0.8635 & 287.0 & 0.8536 & 0.4458 & {[}0.0045, 0.0051{]} \\
         & 30000 & 0.8417 & 0.8497 & 3045 & 0.8411 & 1.232 & {[}-0.0006, -0.0001{]} \\
        \multirow{-5}{*}{250} & 50000 & \cellcolor[HTML]{EFEFEF}0.8039 & \cellcolor[HTML]{EFEFEF}0.8484 & \cellcolor[HTML]{EFEFEF}3600 & 0.8369 & 1.869 & \colorbox{tablegraybg}{{[}0.0286, 0.0295{]}} \\ \cmidrule(l){2-8} 
    \end{tabular}
\end{table}

\begin{table}[ptbh]
    \centering
    \caption{Comparison of MIP and Greedy, under RAEPC spread on 6-core EU email network with $t_0=1$, $p=0.95$ and $r=0.01$: computational time and 95\% confidence interval on the difference $\mu = f(S^h_n) - f(S^\ast_n)$ using McNemar's procedure ($n^{\prime\prime}=2000000$). A positive difference indicates that the greedy solution performed better. Note that columns $\bm{\underline{z}^{\ast}_n}$, $\bm{\overline{z}^{\ast}_n}$, $\bm{z^h_n}$ and both wall times are not from the McNemar's procedure.}
    \label{tab:raepc_compare_mip_and_greedy_t1_p95_r1_appendix}
    \begin{tabular}{@{}cccccccc@{}}
        \toprule
        \multicolumn{2}{c}{\textbf{}} & \multicolumn{3}{|c|}{\textbf{MIP}} & \multicolumn{2}{c|}{\textbf{Greedy}} & $\bm{f(S^h_n) - f(S^\ast_n)}$ \\ \midrule
        $\bm{k}$ & $\bm{n}$ & $\bm{\underline{z}^{\ast}_n}$ & $\bm{\overline{z}^{\ast}_n}$ & \textbf{\begin{tabular}[c]{@{}c@{}}Wall time\\ (sec)\end{tabular}} & $\bm{z^h_n}$ & \textbf{\begin{tabular}[c]{@{}c@{}}Wall time\\ (sec)\end{tabular}} & \textbf{95\% CI} \\ \midrule
         & 1000 & 0.8930 & 0.8990 & 6.972 & 0.8850 & 0.095 & {[}0.0038, 0.0047{]} \\
         & 5000 & 0.8440 & 0.8524 & 1266 & 0.8404 & 0.1276 & {[}0.0003, 0.0011{]} \\
        \multirow{-3}{*}{50} & 10000 & 0.8344 & 0.8427 & 1377 & 0.8329 & 0.2478 & {[}-0.0017, -0.0011{]} \\ \cmidrule(l){2-8} 
         & 1000 & 0.9920 & 1 & 0.2661 & 0.9920 & 0.1224 & {[}0.0271, 0.0281{]} \\
         & 5000 & 0.9608 & 0.9704 & 1477 & 0.9564 & 0.1950 & {[}-0.0020, -0.0013{]} \\
        \multirow{-3}{*}{100} & 10000 & \cellcolor[HTML]{EFEFEF}0.9476 & \cellcolor[HTML]{EFEFEF}0.9613 & \cellcolor[HTML]{EFEFEF}3600 & 0.9454 & 0.3036 & \colorbox{tablegraybg}{{[}-0.0031, -0.0025{]}} \\ \cmidrule(l){2-8} 
    \end{tabular}
\end{table}

\section{Optimal placement in a wheel network}
\label{sec:fn_wheelnetworkplacement_appendix}
We use the term ``node X'' (or just X) to mean the node where detector X has been placed. We compute the probabilities of detecting a virus when detectors are placed in the four configurations shown in Figure~\ref{fig:detector_placements}. Let $A_{\text{adj}}$ be the event that the virus is detected with adjacent detectors, $A_{\text{do}}$ be the event that the virus is detected with diametrically-opposite detectors, $A_{\text{ce}}$ be the event that the virus is detected with center-edge detectors, and $A_{\text{cc}}$ be the event that the virus is detected with center-center detectors, under the network assumptions outlined in Section~\ref{sec:fn_wheelnetworkplacement}. Then we have:
\begin{align*}\label{eq:detectoradjacent_appendix}
    \mathds{P}(A_{\text{adj}}) = & \ \mathds{P}(\text{virus hops from D1 or D2 to a node in }\{D1, D2 \}^C) \cdot (1-r) \\
    & + \mathds{P}(\text{virus hops from D1 to D2 or from D2 to D1}) \cdot (1-r^2) \\
    & + \mathds{P}(\text{virus hops from the center to either D1 or D2})\cdot (1-r) \\
    & + \mathds{P}(\text{virus hops from a neighboring edge node of D1,}\\
    & \quad \text{other than D2, to D1, or from a neighboring edge node} \\
    & \quad \text{of D2, other than D1, to D2)} \cdot (1-r) \\
    = & \ \frac{2}{3}\cdot \frac{2}{v}\cdot (1-r) + \frac{1}{3}\cdot \frac{2}{v}\cdot (1-r^2) + \frac{2}{v-1}\cdot \frac{1}{v}\cdot (1-r) + 2 \cdot \frac{1}{3} \cdot \frac{1}{v}\cdot (1-r) \\
    =  & \ \frac{2}{3}\cdot \frac{2}{v}\cdot (1-r) + \frac{1}{3}\cdot \frac{2}{v}\cdot (1-r^2) \\
    & + \frac{2}{v-1}\cdot \frac{1}{v}\cdot (1-r) + 2 \cdot \frac{1}{3} \cdot \frac{1}{v}\cdot (1-r) \\
    = & \ \frac{(1-r)}{3v(v-1)} \left[6v+2v-2+2rv-2r \right] \\
    = & \ \frac{2(1-r)}{3v(v-1)} \left[4v-1+rv-r \right],
\end{align*}

\begin{equation*}\label{eq:detectordiaopposite_appendix}
\begin{aligned}
    \mathds{P}(A_{\text{do}}) =  & \ \mathds{P}(\text{virus is at D1 or D3 at  } t=0)\cdot (1-r) \\
    & + \mathds{P}(\text{virus hops from the center     to either D1 or D3})\cdot (1-r) \\
    & + \mathds{P}(\text{virus hops from a  neighboring edge node of D1 to D1,} \\
    & \quad \text{or from a neighboring edge     node of D3 to D3})\cdot (1-r) \\
    = & \ \frac{2}{v}\cdot (1-r) +     \frac{2}{v-1}\cdot \frac{1}{v}\cdot (1-r) + 2     \cdot \frac{1}{3}\cdot \frac{2}{v}\cdot (1-r) \\
    = & \ \frac{2(1-r)}{3v(v-1)} \left(5v-2 \right),
\end{aligned}
\end{equation*}

\begin{equation*}\label{eq:detectorcenteredge_appendix}
\begin{aligned}
    \mathds{P}(A_{\text{ce}}) =  & \ \mathds{P}(\text{virus hopping from the   center (D5) to an edge other than D1})\cdot   (1-r) \\
    & + \mathds{P}(\text{virus hopping from the     center to D1})\cdot (1-r^2) \\
    & + \mathds{P}(\text{virus hopping from D1 to an edge node})\cdot (1-r) \\
    & + \mathds{P}(\text{virus hopping from D1 to the center})\cdot (1-r^2) \\
    & + \mathds{P}(\text{virus hopping from a   neighboring edge node of D1} \\
    & \quad \text{ to either D1 or D5}) \cdot (1-r) \\
    & + \mathds{P}(\text{virus hopping from }     \text{(neighbors of D1)}^C - \{D1, D5\} \text{ to D5})\cdot (1-r) \\
    = & \ \frac{v-2}{v-1}\cdot \frac{1}{v}\cdot    (1-r) + \frac{1}{v-1}\cdot \frac{1}{v}\cdot   (1-r^2) + \frac{2}{3}\cdot \frac{1}{v}\cdot     (1-r) \\
    & + \frac{1}{3}\cdot \frac{1}{v}\cdot (1-r^2) + \frac{2}{3}\cdot \frac{2}{v}\cdot (1-r) +     \frac{1}{3}\cdot \frac{v-4}{v}\cdot (1-r) \\
    = & \frac{1-r}{v}\left[\frac{v^2+rv+3v-v-r-3+3 v+3r-3}{3(v-1)} \right] \\
    = & \ \frac{(1-r)}{3v(v-1)}\left[v^2+5v+rv+2r-6  \right],
\end{aligned}
\end{equation*}
and
\begin{equation*}\label{eq:detectorcentercenter_appendix}
\begin{aligned}
    \mathds{P}(A_{\text{cc}}) =  & \ \mathds{P}(\text{virus is at the center   at } t=0)\cdot (1-r^2) \\
    & + \mathds{P}(\text{virus hops from an edge    node to the center})\cdot (1-r^2) \\
    =  & \ \frac{1}{v}(1-r^2) + \frac{1}{3}\cdot     \frac{v-1}{v}\cdot (1-r^2) \\
    =  & \ (1-r^2) \left[\frac{v+2}{3v} \right] \\
    = & \ \frac{(1-r)(1+r)(v+2)}{3v}.
\end{aligned}
\end{equation*}

\subsection{Comparison}
\label{sec:wheelnetworkcomparison_appendix}
\noindent\textbf{Diametrically-opposite vs.\ adjacent detector placements:} 
We now derive a necessary and sufficient condition for the diametrically-opposite placement of detectors to be at least as good as the 
adjacent placement of detectors in the wheel network. 
Using Table~\ref{tab:detectionprobability},
we compare the corresponding detection probabilities, deriving the following condition:
\begin{equation}\label{eq:comparisondiaoppositeandadjacent}
    \begin{aligned}
    & & 5v-2 \geq\ & 4v-1+rv-r \\
    \Leftrightarrow & & v-1+r-rv \geq\ & 0 \\
    \Leftrightarrow & & (v-1)(1-r) \geq\ & 0 \\
    \Rightarrow & & v \geq\ & 1\qquad (\because 0\leq r < 1). \\
    \end{aligned}
\end{equation}
\begin{result}
	From (\ref{eq:comparisondiaoppositeandadjacent}) we see that placing detectors on diametrically-opposite edge nodes is always better than placing them on adjacent-edge nodes.
\end{result}

\noindent\textbf{Diametrically-opposite vs.\ center-edge detector placements:} We next derive conditions to determine when the diametrically-opposite placement is at least as good as the center-edge placement of detectors. We again use 
Table~\ref{tab:detectionprobability} to make the comparison:
\begin{equation}\label{eq:comparisondiaoppositecenteredge1}
    \begin{aligned}
    & & 2(5v-2) \geq\ & v^2+5v+rv+2r-6 \\
    \Leftrightarrow & & v^2-5v+rv+2r-2 \leq\ & 0 \\
    \Leftrightarrow & & v^2+ (r-5)v+ (2r-2) \leq\ & 0. \\
    \end{aligned}
\end{equation}
When $r=0$, the left-hand side of the final expression in (\ref{eq:comparisondiaoppositecenteredge1}) reduces to $v^2 -5v -2$. The roots of $v^2 -5v -2$ are $\frac{5 \pm \sqrt{33}}{2}$. The larger root is smaller than 5.5, since $\frac{5 + \sqrt{33}}{2} < \frac{5 + \sqrt{36}}{2}=5.5$. Because the coefficient of the quadratic term is positive, $v^2 -5v -2>0$ for all $v\geq 6$. Since the left-hand side of the final expression in (\ref{eq:comparisondiaoppositecenteredge1}) is increasing in $r$, this holds for $r \in [0, 1)$.
\begin{result}
	For $v\geq 6$ the center-edge placement is better than placing the detectors on diametrically-opposite edge nodes. For $v=5$, straightforward algebra implies that if $r>\frac{2}{7}$, then the center-edge placement is better, and if $r<\frac{2}{7}$, then the diametrically-opposite placement is better.
\end{result}

\noindent\textbf{Center-edge vs.\ center-center detector placements:} 
We now derive conditions to determine when the center-center placement is at least as good as the center-edge placement. We use expressions from
Table~\ref{tab:detectionprobability} to make the comparison:
\begin{equation}
    \begin{aligned}\label{eq:comparisoncentercenter_centeredge1}
    & & \frac{v^2+5v+rv+2r-6}{v-1} \leq\ & (1+r)(v+2) \\
    \Leftrightarrow & & v^2+5v+rv+2r-6 \leq\ & v^2+2v+rv^2+2rv-v-2-rv-2r \\
    \Leftrightarrow & & 5v+rv+2r-6 \leq\ & v+rv^2+rv-2r-2 \\
    \Leftrightarrow & & rv^2-4v-4r+4 \geq\ & 0. \\
    \end{aligned}
\end{equation}
The roots of the expression $rv^2-4v-4r+4$, viewed as a function of $v$, are $\frac{2}{r} \pm 2\sqrt{1 - \frac{1}{r} + \frac{1}{r^2}}$. This means that when $r>0$, the left-hand side of the final expression in (\ref{eq:comparisoncentercenter_centeredge1}) is non-negative for $v \ge \frac{2}{r} + 2\sqrt{1 - \frac{1}{r} + \frac{1}{r^2}}$. Therefore, when $r=\frac{16}{21}$, then the left-hand side of the final expression in (\ref{eq:comparisoncentercenter_centeredge1}) is non-negative for $v \ge 5$. It can be shown that $\frac{2}{r} + 2\sqrt{1 - \frac{1}{r} + \frac{1}{r^2}}$ is decreasing in $r$ for $r \in [0, 1)$, so the claim holds for $r \ge \frac{16}{21}$.
\begin{result}
	For $r \ge \frac{16}{21} \approx 0.76$ and $v \ge 5$, the center-center placement is as least as good as the center-edge placement.
\end{result}

\pagebreak
\bibliographystyle{unsrtnat}
 \newcommand{\noop}[1]{}

\end{document}